\documentclass [11pt,twoside,a4paper]{article}
\usepackage{amsfonts}
\usepackage{amsthm}
\usepackage{amsmath}
\usepackage{amstext}
\usepackage{amssymb}
\usepackage{mathrsfs}
\usepackage{amscd}
\usepackage{xypic}
\usepackage{epsf}              
\usepackage{graphicx}          
\usepackage{fancybox}          
\usepackage{color}             
\usepackage{fancyhdr}
\usepackage[hang,footnotesize]{caption2}  
\usepackage{mathbbold} 

\def\mathcal{\mathscr}
\newfont{\aaa}{cmb10 at 19pt}
\newfont{\bbb}{cmb10 at 11pt}
\newtheorem{lem}{Lemma}[]
\newtheorem{thm}{Theorem}
\newtheorem{cor}{Corollary}[]
\newtheorem{rem}{Remark}[]
\newtheorem{pro}{Proposition}[]
\def\dsum{\displaystyle\sum}

\def\dsup{\displaystyle\sup}

\newcommand{\beq}{\begin{equation}}
\newcommand{\eeq}{\end{equation}}
\newcommand{\bey}{\begin{eqnarray}}
\newcommand{\eey}{\end{eqnarray}}
\newcommand{\beyy}{\begin{eqnarray*}}
\newcommand{\eeyy}{\end{eqnarray*}}

\makeatletter
\def\@evenhead{
\vbox{\hbox to \textwidth {}{\hspace{0mm}{\footnotesize
\thepage}}{\hspace{5.5cm} {\footnotesize{Mu-Fa CHEN, Ling-Di WANG, Yu-Hui ZHANG}}} \protect\vspace{1truemm}\relax \hrule depth0pt
height0.15truemm width\textwidth}}
\def\@evenfoot{}
\def\@oddhead{\vbox{\hbox to \textwidth
{{\hspace{0cm}{\footnotesize Mixed Principal Eigenvalues in Dimension One}\hfill{\footnotesize
\thepage}}\hspace{0mm}}{} \protect\vspace{1truemm}\relax\hrule
depth0pt height0.15truemm width\textwidth}}
\def\@oddfoot{}
\makeatother

\renewcommand{\thefootnote}{\fnsymbol{footnote}}
\def\scr{\mathscr}
\def\dsup{\sup}
\def\dsum{\sum}
\def\d{\text{\rm d}}
\def\supp{\text{\scriptsize\rm supp}}

\newcommand{\rf}[2]{[\ref{#1}; #2]}

\definecolor{1}{rgb}{1,0,0}
 \definecolor{2}{rgb}{0,1,0}
 \definecolor{3}{rgb}{0,0,1}

 \definecolor{1red}   {rgb}{0.70,0.00,0.00}
 \definecolor{1green} {rgb}{0.00,0.50,0.00}
 \definecolor{1blue}  {rgb}{0.00,0.00,0.70}
 \definecolor{1brown} {rgb}{0.60,0.30,0.00}
 \definecolor{1pink}  {rgb}{0.40,0.00,0.40}
 \definecolor{1cyan}  {rgb}{0.00,0.30,0.30}
 \definecolor{1orange}{rgb}{0.80,0.40,0.00}
 \definecolor{1gray}  {rgb}{0.30,0.30,0.30}

\begin{document}

\thispagestyle{empty} \thispagestyle{fancy} {
\fancyhead[RO,LE]{\footnotesize Front. Math. China 2012,\\[3mm]} \fancyfoot[CE,CO]{}}
\renewcommand{\headrulewidth}{0pt}


\qquad\\[8mm]

\noindent{\aaa{Mixed Principal Eigenvalues in Dimension One}\\[1mm]

\noindent{\bbb Mu-Fa CHEN,\quad Ling-Di WANG,\quad Yu-Hui ZHANG}\\[-1mm]

\noindent\footnotesize{School of Mathematics Sciences, Beijing Normal University, Laboratory of Mathematics and Complex Systems (Beijing Normal University),  Ministry of Education, Beijing 100875, China}\\[6mm]

\vskip-2mm \noindent{\footnotesize$\copyright$ Higher Education
Press and Springer-Verlag Berlin Heidelberg 2012} \vskip 4mm

\normalsize\noindent{\bbb Abstract}\quad This is one of a series of papers exploring the stability speed of one-dimensional stochastic processes. The present paper emphasizes on the principal eigenvalues of elliptic operators.
 The eigenvalue is just the best constant in the $L^{2}$-Poincar\'e inequality and describes the decay rate of the corresponding diffusion process.
We present some variational formulas for the mixed principal eigenvalues of the operators.
As applications of these formulas, we obtain case by case explicit estimates, a criterion for positivity, and an approximating procedure for the eigenvalue.\vspace{0.3cm}

\footnotetext{Received March 6, 2012; accepted June 8,
2012\\
\hspace*{5.8mm}Corresponding author: Ling-Di Wang, E-mail:
wanglingdi@mail.bnu.edu.cn}

\noindent{\bbb Keywords}\quad Eigenvalue, variational formula, explicit estimates, positivity criterion, approximating procedure\\
{\bbb MSC}\quad 60J60, 34L15\\[0.4cm]
\allowdisplaybreaks[4]
\renewcommand{\thefootnote}{\arabic{footnote}}

\noindent{\bbb{1\quad Introduction}}\\[0.1cm]
This paper is a continuation of \cite{r2} in which the stability speed was carefully studied in the discrete situation (birth--death processes) and partially in the continuous one
(diffusions). For a large part of the study, the description of the problem is equivalent to the Poincar\'e-type inequalities or the principal eigenvalue. On the last two topics,
there are a great number of publications (cf.  \cite{r6, rn1} and
references therein for the background and motivation of the study on these topics). However, to save the space here,  most of the references are not repeated in this paper. Consider a finite interval $(0, D)$ for a moment. We are  interested in some typical Sturm-Liouville eigenvalue problems. According to the Dirichlet
(denoted by code ``D'') and Neumann (denoted by code ``N'')
boundaries at the left- or right-endpoint, we have four cases
of boundary condition:
DD, ND, DN and NN. In the diffusion context, the DD- and NN-cases are largely handled in [\ref{r4}\,--\,\ref{r2}] and \cite{rn3,r1} . The present paper is mainly devoted to the ND- and DN-cases.
As will be seen in the next section, the classification for the boundaries is
also meaningful when $D=\infty$.

 The paper is organized as follows. In the next section, we focus on the ND-case. First, we introduce several variational formulas for the eigenvalue.
As a consequence, we obtain the basic estimates, a criterion for positivity, an approximating procedure, and improved estimates for the eigenvalue.
 As far as we know, most of these results, except Theorem \ref{t11.2.1}, have not yet appeared in the literatures. The proofs of them are sketched in Section 3.
  From \rf{r2}{Section 10}, we know that the DN-case and the ND-case are dual to each other. Thus, as a dual to the ND-case, it is natural to study the DN-case. To which Section 4 is devoted, partial results come from the duality but some of them are not, need direct proofs. The main extension to the earlier study is that here we do not assume the uniqueness of the processes, instead of which we adopt the maximal extension of the Dirichlet form or the maximal process. Finally, some supplement to \cite{r3, r5, r1} in the NN-case (i.e., the ergodic case.) is presented in Section 5. The complete proofs of the results presented in this paper are quite technical and long. However, a large part of them are parallel to \cite{r2} and so we omit mostly the ``translation'' from the discrete situation to the continuous one. Instead, we emphasis on the difference between them (Lemmas \ref{L1.1.0}--\ref{L1.1.4}, for instance), and illustrate a little of the translation for the reader's reference. We may leave the details to our homepage or publish them elsewhere.

  The basic estimates are also studied in \cite{wj} in terms of $H$-transform.

Some examples of the study are illustrated in \rf{rn1}{Section 5}. The most powerful application of the improved estimates presented in the paper is given by
\cite{rn2} where the lower and upper bounds are quite close or almost coincide with each other.

Here we discuss briefly about the problem on the whole line.
We consider the ND-case only.
First, one may regard the whole line $\mathbb R$ as a limit of $[M, \infty)$ as $M$ decreases to $-\infty$. Then the mixed eigenvalue problem on $[M, \infty)$ is known by what we are studying in the paper. Next, one may split $\mathbb R$
into two parts: $(-\infty, 0)$ and $(0, \infty)$.
The case with ND-boundaries on $(0, \infty)$ is studied in Sections 2 and 3. Besides, the case with ND-boundaries on $(-\infty, 0)$ is simply a reverse of the DN-case on $(0, \infty)$. Therefore, the behavior of the original operator
on the whole line should be clear. However, there is an interesting point here. On $(0, \infty)$, we use the minimal Dirichlet form but on $(-\infty, 0)$ we adopt the maximal one. Thus, the domain of the original Dirichlet form on the whole line may be neither
the maximal nor the minimal one. Therefore, it is essentially different from
DD- or NN-cases on the whole line we have studied in \cite{r2},
\cite{rn1} and \cite{rn3}.

To conclude this section, we mention that in a more general context, for the Poincar\'e-type inequalities, the DN-case
was completed earlier (cf. \rf{r6}{Chapter 6}), the basic estimates for the ND-case in the discrete situation was given by \rf{r2}{Theorem 8.5} from which one can write down easily the continuous version.
\noindent \\[4mm]

\noindent{\bbb 2\quad The ND-case}\\[0.1cm]
Define
$$\aligned
{\scr C}[0, D]&=\{f:f\;\text{\rm is continuous on}\;[0, D]\}\quad \text{and}\\
 {\scr C}^k(0, D)&=\{f: f \text{ has continuous derivatives of order $k$ on $(0, D)$}\},\qquad k\geqslant 1.
 \endaligned$$
Here and in what follows, when $D=\infty$, the notation
${\scr C}[0, D]$ simply means ${\scr C}[0, D)$. The convention
should be clear in other cases but we will not mention time by time.
Let
$$L=a(x)\frac{\d^{2}}{\d x^{2}}+b(x)\frac{\d}{\d x}$$ be an elliptic operator on an interval $(0,D)\,(D\leqslant\infty)$. Set
 $$C(x)=\int_{0}^{x} \frac{b(u)}{a(u)}\d u .$$
Throughout this paper, we need the following hypothesis
(which is trivial in the discrete situation):
\begin{align}
 &\text{The functions $a$, $b$ are Borel measurable on [0,D] and $a$ is positive on $(0,D)$,}\nonumber\\
  &\text{$b/a$ and $e^C/a$ are locally integrable on $[0,D]$.}
\label{f1.0}
\end{align}
Note that for continuous $a$ and $b$, the hypothesis \eqref{f1.0} is
reduced to the condition $a>0$ only.
In this section, we consider the ND-boundaries only.
More precisely, as usual, the Dirichlet boundary condition at $D$ means that $g(D)=0$ when $D<\infty$. When $D=\infty$, it is natural to take ``$\lim_{x\to\infty} g(x)=0$'' as a boundary condition. However, this is not pre-assumed but proved later
(cf.  Lemma \ref{L1.1.4} below). Therefore, the code ``ND'' is still
meaningful even if $D=\infty$.

Throughout this section, we work on the following mixed
principal eigenvalue:
\begin{equation}
\lambda_{0}=\inf\big\{D(f):\mu\big(f^{2}\big)=1,\,f\in{\scr C}_K[0, D],\, f(D)=0\text{ if }D<\infty\big\}, \label{f1.1}\end{equation}
where $\mu(f):=\int_{0}^{D} f \d \mu,$
$$\aligned&{\scr C}_K [0,D]=\big\{f:f\in {\scr C}^1(0, D)\cap {\scr C}[0, D],\, f\;\text{\rm has compact support}\big\},\\
&D(f)=\int_{0}^{D}a {f'}^2\d\mu,\qquad \mu(\d x)=\frac{e^{C(x)}}{a(x)}\d x.
\endaligned$$
Besides $\mu$, throughout the paper, we often use another measure:
$\nu(\d x)=e^{-C(x)}\d x$.
When $D<\infty$, $\lambda_0$ coincides with the minimal solution $\lambda$ to the following eigenequation
$$ L f =- \lambda f, \qquad f'(0)=0, \qquad\text{and\quad $f(D)=0$ if $D<\infty$}.$$

To state our results, we need some notation. Define
$$\aligned
&I(f)(x)=-\frac{e^{-C(x)}}{f'(x)}\int_{0}^{x}f\,\d\mu \qquad \mbox{(single integral form)},\\
&I\!I(f)(x)=\frac{1}{f(x)}\!\int_{(x, D)\cap \supp  (f)}\hskip -0.2cm\nu(\d s)\!\int_{0}^{s}\!\!f\d \mu,\;x\in\text{supp}(f)\;\,\mbox{(double integral form)},\\
&R(h)(x)=-\big(ah^{2}+bh+ah'\big)(x)\qquad\mbox{(differential form)}.\endaligned$$
 The domains of the three operators defined above are, respectively, as follows.
$$\aligned
&\mathscr{F}_{I}=\big\{f: f\in {\scr C}^1(0, D)\cap {\scr C}[0, D],\, f|_{(0, D)}>0,\text{ and }f'|_{(0,D)}<0\big\},\\
&\mathscr{F}_{I\!I} =\{f: f\in {\scr C}[0,D]\text{ and } f|_{(0, D)}>0\},\\
&\mathscr{H}=\big\{h: h\in {\scr C}^1(0, D)\cap {\scr C}[0, D],\, h(0)=0,\,h|_{(0, D)}<0\text{ if }\nu(0, D)<\infty,\\
&\qquad\qquad\, \text{and } h|_{(0, D)}\leqslant 0\text{ if }\nu(0, D)=\infty\big\},\quad \mbox{$\nu(\alpha, \beta):=\int_{\alpha}^\beta \d\nu$}.
\endaligned$$
These sets are used for the lower estimates of $\lambda_0$. For the upper bounds, some modifications are needed to avoid the non-integrability  problem, as shown below.
$$\aligned
&\mathscr{{\widetilde F}}_{I}\!=\big\{f\!:\! f\!\in\! {\scr C}^1\!(x_{0}, x_{1}) \cap {\scr C}[x_{0},x_{1}], f'|_{(x_{0},x_{1})}\!<0\, \text{for some}\, x_0,x_1\!\in\![0,D)\text{ with }\\&\text{\hskip 4em}\,x_0\!<\!x_1,
\text{and } f=f(\cdot\vee x_{0})
\mathbbold{1}_{[0, x_{1})}\big\},\\
&\mathscr{{\widetilde F}}_{I\!I}\!=\big\{f: \exists x_{0}\in(0, D)\text{ such that } f=f \mathbbold{1}_{[0, x_{0})}\text{ and } f\in  {\scr C}[0,x_{0}]\big\},\\
&\mathscr{\widetilde{H}}\!=\!\Big\{h\!:\! \exists x_{0}\!\in\!(0,D)\text{ such that } h\!\in\!{\scr C}^1(0,x_{0})\cap{\scr C}[0,x_{0}], h|_{(0, x_0)}\!\!<\!0,\,h|_{[x_0,D]}\!=\!0,\\
&\text{\hskip 4em}\text{ and } h(0)=0,\,\sup_{(0, x_0)}\big(ah^{2}+bh+ah'\big)<0\Big\}.
\endaligned$$
Here and in what follows, we adopt the usual convention $1/0=\infty$. The superscript ``$\widetilde{\quad}$'' means modified. In the formulas of Theorem \ref{th1}  below, ``$\sup\,\inf$'' are
used for lower bounds of $\lambda_0$, each test function $f$
produces a lower bound $\inf_x I(f)(x)^{-1}$, and so this part
is called variational formula for the lower estimate of $\lambda_0$. Dually, the
``$\inf\,\sup$'' are used for upper estimates of $\lambda_0$. Among them, the ones expressed by the operator $R$ are easiest
to compute in practice, and the ones expressed by $I\!I$
are hardest to compute but provide better estimates.
Because of ``$\inf\,\sup$'', a localizing procedure is used
for the test function to avoid $I(f)\equiv \infty$ for instance,
which is removed out automatically for the ``$\sup\,\inf$'' part. Each part of  Theorem \ref{th1} below plays a role in our study.
Parts (1) and (2) are applied to Theorems \ref{t11.2.1} and \ref{t11.2.2}, respectively. Part (3) is a
comparison with Proposition \ref{t11.1.3}, which is then used for a dual form of Theorem \ref{1t1.1}\,(3).

\begin{thm}\label{th1} Under hypothesis \eqref{f1.0}, the following variational formulas hold for $\lambda_{0}$ defined by \eqref{f1.1}.

$(1)$ Single integral forms:
$$
\aligned
\inf\limits_{f\in\mathscr{{\widetilde F}}_{I}}\,\dsup_{x\in(0, D)}I(f)(x)^{-1}=\lambda_{0}=\sup\limits_{f\in\mathscr{F}_{I}}\,\inf\limits_{x\in(0,D)}I(f)(x)^{-1},
\endaligned$$

$(2)$ Double integral forms:
$$
\aligned
\inf\limits_{f\in\mathscr{{\widetilde F}}_{I\!I}}\,\dsup_{x\in
\supp(f)}I\!I(f)(x)^{-1}=\lambda_{0}=\sup\limits_{f\in\mathscr{F}_{I\!I}}\,\inf\limits_{x\in(0, D)}I\!I(f)(x)^{-1}.
\endaligned$$
Moreover, if $a,\, b\in {\scr C}[0,D]$, then we have additionally a

$(3)$  differential form:
$$
\aligned
 \inf\limits_{h\in\mathscr{\widetilde{H}}}\,\dsup_{x\in(0,
D)}R(h)(x)=\lambda_{0}=\sup\limits_{h\in\mathscr{H}}\,\inf\limits_{x\in(0, D)}R(h)(x).
\endaligned$$ Furthermore, the supremum on the right-hand side of the above three formulas  can be attained.
\end{thm}

The next result, similar to the discrete case, either extends the domain of
$\lambda_{0}$, or adds some additional sets of test functions for operators $I$ and $I\!I$, respectively. Besides, as an application of the lower variational formula (Theorem \ref{th1}\,(2)), we obtain the vanishing property of the eigenfunction (Lemma \ref{L1.1.4})
which leads to the crucial part (1) of the proposition below. The vanishing property is the meaning of the Dirichlet boundary at $D=\infty$ as we expected. A more common description of $\lambda_{0}$ is given by Lemma \ref{new} below.

\begin{pro}\label{p11.1.1}
Let hypothesis \eqref{f1.0} hold. Then

$(1)$ we have
$$\aligned
\lambda_{0}&=\inf\Big\{D(f):\mu\big(f^{2}\big)=1,\,f\in {\scr C}^1(0, D)\cap {\scr C}[0,D]\text{ and }f(D)=0\Big\}\\
&=:\tilde{\lambda}_{0},
\endaligned$$
where $f(D)=\lim_{x\to D}f(x)$ in the case of $D=\infty.$

$(2)$ Moreover, we have
\begin{eqnarray}
&&\inf\limits_{f\in\mathscr{{\widetilde F}'}_{I}}\,\dsup_{x\in(0,
D)}I(f)(x)^{-1}=\lambda_{0}=\sup\limits_{f\in\mathscr{F}_{I}}\,\inf_{x\in(0, D)}I\!I(f)(x)^{-1},\label{fn1}\\
&&\inf\limits_{f\in\mathscr{{\widetilde F}}_{I\!I}\bigcup\mathscr{{\widetilde F}}_{I\!I}'}
\,\dsup_{x\in \supp(f)}I\!I(f)(x)^{-1}=\lambda_{0}=\inf\limits_{f\in\mathscr{{\widetilde F}}_{I}}\,\dsup_{x\in
\supp(f)}I\!I(f)(x)^{-1},\label{fn9}
\end{eqnarray}
where $$\aligned
&\mathscr{{\widetilde F}}_{I}'= \big\{f: \exists x_{0}\in(0, D)\text{ such that }f=f \mathbbold{1}_{[0, x_{0})},\, f\in {\scr C}^1(0, x_{0})\\
&\text{\hskip 4em}\cap {\scr C}[0,x_{0}],\text{ and } f'|_{(0, x_{0})}<0\big\},\\
&\mathscr{{\widetilde F}'}_{I\!I}=\big\{f: f>0,\, f\in {\scr C}[0,D],\text{ and } fI\!I(f)\in L^{2}(\mu)\big\}.
\endaligned$$
Besides, the supremum over $\{f\in\mathscr{F}_{I}\}$ in \eqref{fn1} can be attained.
\end{pro}
The operator $\overline R$ defined below was first introduced in \rf{r1}{Theorem 2.1} based on a probabilistic (coupling) technique. Different from $R$, it is a ``bridge'' in proving the duality of the ND- and DN-cases. It also leads to a different variational formula for $\lambda_0$ as follows.
\begin{pro}\label{t11.1.3} Suppose that $a, b\in {\scr C}^1 (0, D)\cap {\scr C} [0, D]$ and $a>0$ on $(0,D)$. Set
$$\overline{\mathscr{H}}=\big\{h:  h(0)=0,\,h\in {\scr C}^2(0,D)\cap\mathscr{C}[0,D],\text{ and } h|_{(0, D)}<0\big\}$$
and define
$${\overline R}(h)(x)=-\frac{(ah'+b h)'(x)}{h(x)}. $$
Then

 $(1)$  we have
 $\sup_{h\in\overline{\mathscr{H}}}\inf_{x\in(0,  D)}{\overline R}(h)(x)\geqslant\lambda_{0}$ and the equality sign holds once  $\mu(0, D)=\infty$.

 $(2)$ In general, we have
 \begin{equation}\label{fn2}\lambda_{0}=\sup\limits_{h\in\mathscr{H}_{\ast}}\inf\limits_{x\in(0,D)}\overline {R}(h)(x),\end{equation}
where
$$\mathscr{H}_{\ast}=\big\{h\in {\scr C}^2(0, D)\cap {\scr C}[0,D]: h(0)=0, \text{ and } h<0,\,h'<-{a}^{-1}{b}h\text{ on }(0,D)\big\}.$$
Moreover, the supremum in \eqref{fn2} can be attained.\end{pro}

\begin{rem}\label{rem1}{\rm (Comparison of $R$ and $\overline R$ )} With
$h={g'}/{g}$, we have
$$-\frac{L g}{g}=-(ah^{2}+bh+ah')=R(h).$$
Next, with $h=g'$, we have
$$-\frac{(L g)'}{g'}=-\frac{(ah'+b h)'}{h}=\overline R(h).$$
\end{rem}

As an application of Theorem \ref{th1}\,(1) to the test function
$\nu(x, D)^{\gamma}$ with $\gamma=1/2$ or $ 1$, we obtain the basic estimates and furthermore a criterion as follows.

\begin{thm}\label{t11.2.1} $(\text{\rm Criterion and basic estimates})$
Let hypothesis \eqref{f1.0} hold.
Then
$\lambda_{0}>0$ iff
$$\delta:=\sup_{x\in(0, D)} \mu (0, x)\,\nu (x, D)<\infty,
\qquad \mbox{$\mu(\alpha,\,\beta):=\int_{\alpha}^{\beta}\d\mu.$}$$
More precisely, we have
$$(4\delta)^{-1}\leqslant\lambda_{0}\leqslant\delta^{-1}.$$
In particular, when $D=\infty$, we have
$\lambda_{0}=0$ if $\nu(0, D)=\infty$, and $\lambda_{0}>0$ if
\newline
$\int_0^\infty\mu(0, x)\,\nu(\d x)<\infty$.
\end{thm}

The next result is an application of Theorem \ref{th1}\,(2), repeated with $f=f_n$, starting from the initial $f_1$ , the test
function just mentioned above Theorem \ref{t11.2.1}. The result provides us a way to improve the basic
estimates step by step. In view of the last criterion,
for any improvement, one may assume that $\delta<\infty$.

\begin{thm}\label{t11.2.2}$(\text{\rm Approximating procedure})$ Let hypothesis \eqref{f1.0} hold and assume that $\delta<\infty$.
Set $\varphi(x)=\nu(x,D)$.

$(1)$ Let $f_{1}=\sqrt{\varphi}$, $f_{n}=f_{n-1}I\!I(f_{n-1})$, and $\delta_{n}=\sup_{x\in(0, D)}I\!I(f_{n})(x)$. Then $\delta_{n}$ is decreasing in $n$ and $$ \lambda_{0}\geqslant\delta_{n}^{-1}\geqslant({4\delta})^{-1},\qquad n\geqslant1.$$

$(2)$ For fixed $x_{0},\ x_{1}\in [0, D)$ with $x_{0}<x_{1}$, define
$$f_{1}^{x_{0}, x_{1}}\!=\!
\nu (\cdot\vee x_{0},\, x_1)\, \mathbbold{1}_{[0, x_{1})}, \qquad f_{n}^{x_{0},x_{1}}\!=\!\big(f_{n-1}^{x_{0}, x_{1}}I\!I\big(f_{n-1}^{x_{0},x_{1}}\big)\big)(\cdot\vee x_0)\,\mathbbold{1}_{[0, x_{1})},\quad n\geqslant1,$$
and let $\delta_{n}'=\sup_{x_{0}, x_{1}:x_{0}<x_{1}}\,\inf_{x< x_{1}}I\!I({f_{n}}^{x_{0}, x_{1}})(x)$. Then
$\delta^{-1}\geqslant
{\delta_{n}'}^{-1}\geqslant\lambda_{0}$ for $n\geqslant1.$
\medskip

$(3)$ Define
$$\bar{\delta}_{n}=\sup\limits_{x_{0},x_{1}:x_{0}<x_{1}}\frac{\|{f_{n}}^{x_{0}, x_{1}}\|}{D({f_{n}}^{x_{0}, x_{1}})},\qquad n\geqslant 1. $$
Then $\bar{\delta}_{n}^{-1}\geqslant \lambda_{0},\ \bar{\delta}_{n+1}\geqslant\delta_{n}'\, (n\geqslant 1)$,  and $\bar{\delta}_{1}=\delta_{1}'$.
\end{thm}

The next result comes from the first step of the above approximation.

\begin{cor}\label{cor1}$(\text{\rm Improved estimates})$ Let hypothesis \eqref{f1.0} hold. For $\lambda_{0}$, we have
$$\aligned
\delta^{-1}\geqslant{\delta'_{1}}^{-1}\geqslant\lambda_{0}\geqslant\delta_{1}^{-1}\geqslant(4\delta)^{-1},
\endaligned$$
where
\begin{eqnarray}
\delta_{1}\hskip-0.6cm&&=\sup\limits_{x\in(0,D)}
\frac{1}{\sqrt{\varphi(x)}}\int_{0}^{D}\sqrt{\varphi}\,
\varphi(\cdot\vee x)\,\d\mu\nonumber\\&&=\sup\limits_{x\in(0,D)}
\bigg(\sqrt{\varphi(x)}\,\int_{0}^{x}\sqrt{\varphi}\,\d\mu
+\frac{1}{\sqrt{\varphi(x)}}\int_{x}^{D}\varphi^{{3}/{2}}\d\mu\bigg),\end{eqnarray}
\begin{equation}\label{f1.3}
\delta_{1}'=\sup\limits_{x\in(0,D)}\bigg(\mu(0,x)\,\varphi(x)
+\frac{1}{\varphi(x)}\int_{x}^{D}\varphi^2\,\d\mu\bigg)\in[\delta,2\delta].\end{equation}
\end{cor}
\noindent\\[4mm]
\noindent{\bbb 3\quad Partial proofs of the results in Section 2 }\\[0.1cm]
Some preparations are needed to prove our main results.
The first six lemmas below, except Lemma \ref{new}, are mainly devoted to describe the eigenfunction of $\lambda_{0}$. The Lemmas are essential in our
study.
Note that their proofs are very different from the discrete situation.
The first one below is taken from \rf{zettl}{Theorems 1.2.1 and 2.2.1}.

\begin{lem}\label{L1.1.0}
\begin{itemize}\setlength{\itemsep}{-0.8ex}
\item[$(1)$] Let hypothesis \eqref{f1.0} hold. Then, whenever $g$ and $g'$ are initially not vanished simultaneously, there exists uniquely a non-zero function $g\in {\scr C}^1[0, D]$ such that $g'$ is absolutely continuous on each compact subinterval of $[0, D)$ and the eigenequation
   $L g=-\lambda g$ holds almost everywhere.
\item[$(2)$] Suppose additionally $a$ and $b$ are continuous on $[0,D]$. Then $g\in {\scr C}^2[0, D]$ and the eigenequation holds everywhere on $[0, D]$.
\end{itemize}
 \end{lem}

 In what follows, we call the function $g$ given in part (1) of Lemma \ref{L1.1.0} {\it a.e. eigenfunction of $\lambda$.}
 Remember we need ``a.e.'' only in the case where $g''$ is used.
Of course, we remove ``a.e.'' if the eigenequation holds
everywhere.

The next result enables us to return to a more common description of the eigenvalue.
\begin{lem}\label{new}
Let $\mathscr{A}[\alpha, \beta]$ be the set of all absolutely continuous functions on $[\alpha, \beta]$. Define $$\lambda_{\ast}=\inf\{D(f): f\in\mathscr{A}[0,D],\, \|f\|=1,\;\text{ and } f(D)=0\}.$$ Then $\tilde{\lambda}_{0}=\lambda_{\ast}$.
\end{lem}
{\it Proof}\quad\rm
It is obvious that $\lambda_{\ast}\leqslant\tilde{\lambda}_{0}$. Next, let $g$ be the a.e. eigenfunction of $\lambda_{\ast}$. Then, $g'\in\mathscr{A}[0,D]$ by Lemma \ref{L1.1.0}\,(1). By making inner product with $g$ on the both sides of $Lg=-\lambda_{\ast}g$ with respect to $\mu$, it follows that $$-\big(e^{C}gg'\big)\big|_{0}^{D}+D(g)=\lambda_{\ast}\|g\|^{2}.$$ Since $g'(0)=0$
and $(gg')(D)\leqslant 0$, we have $\lambda_{\ast}\geqslant D(g)/\|g\|^{2}$. Because $g'\in\mathscr{A}[0,D]$, it is clear that $D(g)/\|g\|^{2}\geqslant\tilde{\lambda}_{0}$. We have thus obtained that $$\tilde{\lambda}_{0}\leqslant\lambda_{\ast}\leqslant\tilde{\lambda}_{0},$$
and so $\tilde{\lambda}_{0}=\lambda_{\ast}$. There is a small gap in the proof above since in the case of $D=\infty$, the a.e. eigenfunction $g$ may not belong to $L^{2}(\mu)$ and we have not yet proved that $(gg')(D)\leqslant 0$. However, one may avoid this by a standard approximating procedure,
\footnote{If $D=\infty$ and $\lambda_{\ast}<\tilde{\lambda}_{0}$,
then there would exist $p_{n}<\infty$ such that
$\lambda_{\ast}^{(0,p_n)}<\tilde{\lambda}_0^{(0,p_n)}$
which is a contradiction with what we have just proved.}
using $[0,p_{n}]$ instead of $[0,D)$ with $p_{n}\uparrow D$ provided $D=\infty$:
  $$\aligned
\lim_{n\to\infty}\lambda_{0}^{(0, p_n)}\!\!&
=\lim_{n\to\infty}\!\inf\big\{D\big(f\big): f\in\mathscr{C}[0,p_{n}]\cap{\scr C}^{1}(0,p_{n}),\,\mu\big(f^{2}\big)=1,\,f|_{[p_{n},D]}=0\big\}\\
&=\inf\big\{D(f):\!\mu\big(f^{2}\big)\!=\!1, f\in{\scr C}[0,D]\cap{\scr C}^{1}(0,D), f(D)\!=\!0\text{ if }D<\infty\big\}\\
&=\lambda_{0}^{(0, D)}=\tilde{\lambda}_{0}.
\qquad\Box\endaligned $$

Clearly, because of hypothesis \eqref{f1.0}, we have $\lambda_{0}>0$ once $D<\infty$.
The next result is a simple comparison. For given $\alpha, \beta\,(\alpha< \beta)$, denote by $\lambda_0^{(\alpha, \beta)}$ and $\lambda_1^{(\alpha, \beta)}$ respectively, the principal ND- and NN-eigenvalue (the latter is also called the first nontrivial eigenvalue or the spectral gap in the ergodic case). For simplicity, we use $\downarrow$ (resp. $\downarrow\downarrow,\,\uparrow,\,\uparrow\uparrow)$ to
denote decreasing (resp. strictly decreasing, increasing, strictly increasing).

\begin{lem}\label{L1.1.1}
\begin{itemize}\setlength{\itemsep}{-0.8ex}
\item[$(1)$] For $p, q \in (0, D)$ with $p<q$, we have
$\lambda_0^{(0, p)}> \lambda_0^{(0, q)}$. Furthermore,
$\lambda_0^{(0, p_n)}\downarrow\downarrow \lambda_0^{(0, D)}\,\big(\text{ i.e.}\, \tilde{\lambda}_{0}\big)$ as $p_n\uparrow\uparrow D$.
\item[$(2)$] For $p\in (0, D)$, we have
$\lambda_1^{(0, p)}>\lambda_0^{(0, p)}.$
\end{itemize}
 \end{lem}

{\it Proof}\quad\rm (a) Let $g\,(\ne 0)$ be an a.e. eigenfunction of $\lambda_0^{(0, p)}$. Then $g'(0)=0$, $g(p)=0$ and $Lg=-\lambda_0^{(0, p)} g$ a.e. on $(0, p)$ by Lemma \ref{L1.1.0}\,(1). Moreover
$$\lambda_0^{(0, p)}=\frac{D_{0, p}(g)}{\|g\|_{L^2(0, p;\, \mu)}^2},\qquad D_{\alpha, \beta}(f)=\int_{\alpha}^\beta a {f'}^2 \d\mu.$$
By Lemma \ref{new}, the proof of the first assertion in part (1) will be done once we choose a function $\tilde g\in\mathscr{A}[0,q]$
such that $\tilde g'(0)=0$, $\tilde g(q)=0$, and
\begin{equation}\label{ff3}\frac{D_{0, p}(g)}{\|g\|_{L^2(0, p;\, \mu)}^2}
> \frac{D_{0, q}(\tilde g)}{\|\tilde g\|_{L^2(0, q;\, \mu)}^2}\quad\Big(\geqslant \lambda_0^{(0, q)}\Big).\end{equation}
To do so, without loss of generality, assume that $g|_{(0, p)}>0$ (this is a well-known property as a reverse of the DN-case for finite intervals, cf. \rf{r6}{Theorem 3.7}). Then the required assertion follows
for
$$ {\tilde g}(x)=
(g+\varepsilon)\mathbbold{1}_{[0,p)}(x)+\frac{\varepsilon(x-q)}{(p-q)}\mathbbold{1}_{[p, q]}(x),\qquad x\in[0,q],
$$ once $\varepsilon$ is sufficiently small.
Actually, by simple calculation, we have $$\aligned &D_{0,q}(\tilde g)=D_{0,p}(g)+\frac{\varepsilon^{2}}{(p-q)^{2}} \int_{p}^{q}e^{C(x)}\d x,\\&\|\tilde g\|^{2}_{L^{2}(0,q;\mu)}=\|g\|^{2}_{L^{2}(0,p;\mu)}+\varepsilon\int_{0}^{p}(2g+\varepsilon)\d\mu+\frac{\varepsilon^{2}}{(p-q)^{2}} \int_{p}^{q}(x-q)^{2}\d \mu.\endaligned$$
Thus, \eqref{ff3} holds iff
$$\aligned&\frac{\varepsilon}{(p-q)^{2}}\int_{p}^{q}e^{C}\d x\,\big\|g\big\|_{L^{2}(0,p;\mu)}^{2}\\&\hskip1.5cm<\bigg(\int_{0}^{p}(2g+\varepsilon)\d\mu+\frac{\varepsilon}{(p-q)^{2}}\int_{p}^{q}(x-q)^{2}\d\mu\bigg)D_{0,p}(g).\endaligned$$
Since $\lambda_{0}^{(0,p)}= D_{0,p}(g)\big/\|g\|_{L^{2}(0,p;\mu)}^{2}$, it suffices that
$$\frac{\varepsilon}{(p-q)^{2}}\int_{p}^{q}e^{C}\d x <\lambda_{0}^{(0,p)}\bigg(2\int_{0}^{p}g\d\mu\bigg),$$
which is obvious for sufficiently small $\varepsilon$.

The second assertion in part (1) has just been proved at the end of the last proof.

(b) Part (2) of the Lemma strengthens in the present situation a general result that $\lambda_{1}\geqslant\lambda_{0}$ proved in  \rf{r3}{Proposition 3.2}. Let $g\ne$\,constant be an a.e. eigenfunction of $\lambda_1^{(0, p)}$.
Then $g'(0)=0$, $g'(p)=0$ and $Lg=-\lambda_1^{(0, p)} g$ a.e. on $(0, p)$ by Lemma \ref{L1.1.0}\,(1). Moreover
$$\lambda_1^{(0, p)}=\frac{D_{0, p}(g)}{\text{\rm Var}_{(0, p)}(g)},\qquad \text{\rm Var}_{(\alpha, \beta)}(f)=\int_{\alpha}^\beta f^2\d \mu
-\frac{\mu_{\alpha, \beta}(f)^2}{\mu (\alpha, \beta)}.$$
Without loss of generality, assume that $g$ is strictly increasing (cf. \rf{r2}{Proposition 6.4}).
Then we have
$$\tilde g(x):= g(p)-g(x)>0\qquad \text{on } (0, p).$$
Thus, $\tilde g'(0)=0$, $\tilde g(p)=0$ and moreover
\footnote{Note that $D_{0,p}(\tilde g)=D_{0,p}(g)$ and
 $\text{\rm Var}_{(0,p)}(\tilde g)=\text{\rm Var}_{(0,p)}(g)$.}
$$\lambda_1^{(0, p)}=\frac{D_{0, p}(\tilde g)}{\text{\rm Var}_{(0, p)}(\tilde g)}
=\frac{D_{0, p}(\tilde g)}{\|\tilde g\|_{L^2(0, p;\, \mu)}^2-\mu(\tilde g)^2/\mu(0, p)}
>\frac{D_{0, p}(\tilde g)}{\|\tilde g\|_{L^2(0, p;\, \mu)}^2}
\geqslant \lambda_0^{(0, p)}. \qquad \Box
$$

Before moving further, let us mention a nice expression of $L$:
$$L=\frac{\d}{\d\mu}\frac{\d}{\d\nu}$$
which can be checked by a simple computation. Next, a large part
of the results in the last section is related to the Poisson equation $Lg=-f$, a.e., from which we obtain
\begin{equation}\label{fnew1}\frac{\d}{\d\nu} g(\beta)-\frac{\d}{\d\nu} g(\alpha)=
-\int_\alpha^\beta f \d \mu,\qquad \alpha, \beta \in [0, D],\; \alpha< \beta.\end{equation}
Furthermore, if $g'(\alpha)=0$, then we have \begin{equation}\label{fnew2}g(q)-g(p)=-\int_p^q\nu(\d \beta)\int_\alpha^\beta f\d\mu,\qquad p, q\in[0,D],\; p<q.\end{equation}
Especially, because
$$\frac{\d}{\d\nu} g(0)=e^{C(0)}g'(0)=0,$$
and \eqref{fnew1}, with $f=\lambda_{0}g$, it follows that
\begin{equation}\label{ff1}\frac{\d}{\d\nu} g(s)= -\lambda_0 \int_0^s g \d \mu,\qquad s\in (0, D).\end{equation}

Lemmas \ref{L1.1.2} -- \ref{L1.1.4} given below consist of the basis of the test functions used in the definitions of ${\scr F}_\#$ and ${\scr H}$.
\begin{lem}\label{L1.1.2} Let
$g$ be a non-zero a.e. eigenfunction of $\lambda_{0}>0$. Then $g$ is strictly monotone. \end{lem}

 {\it Proof}\quad\rm Because $\lambda_{0}>0$, $g$ can not be a constant. We need only to prove that $g'\ne 0$ on $(0, D)$.
 Suppose that there is a $p\in (0, D)$ such that $g'(p)=0$.
 Then, by the eigenequation restricted to $(0, p)$, we would have $\lambda_{0}\geqslant \lambda_1^{(0, p)}$, where $\lambda_{1}^{(0,p)}$ is the
 minimal eigenvalue with Neumann boundaries at $0$ and $p$.
 To see this, by (\ref{ff1}),
 we have $\mu_{0,p}(g)=0$ since $g'(0)=0$ and $ g'(p)=0$. From here, it is quite standard
 to prove the required assertion. By making inner product with $g$ on the both sides of
 the eigenequation with respect to $\mu_{0,p}$, it follows that
$$-\big(e^{C}gg'\big)|_{0}^{p}+D_{0,p}(g)=\lambda_{0}\mu_{0,p}\big(g^{2}\big).$$
Again, because of  $g'(0)=g'(p)=0$, we obtain $\lambda_{0}=D_{0,p}(g)/\mu_{0,p}(g^{2}).$
 Hence,
 $$\aligned\lambda_{0}&=\frac{D_{0,p}(g)}{\mu_{0,p}(g^{2})}=\frac{D_{0,p}(g)}{\text{\rm Var}_{(0,p)}(g)}\quad(\text{since } \mu_{0,p}(g)=0)\\
 &\geqslant\inf\bigg\{\frac{D_{0,p}(f)}{\text{\rm Var}_{(0,p)}(f)}: f\in {\scr C}^{1}(0,p)\cap {\scr C}[0,p] ,f\in L^{2}(0,p;\mu),f\neq \text{constant}\bigg\}\\
 &=\lambda_{1}^{(0,p)}.\endaligned$$

 Now, by Lemma \ref{L1.1.1},  we  obtain
 $$\lambda_{0}\geqslant \lambda_1^{(0, p)}>\lambda_0^{(0, p)}
> \lambda_0^{(0, D)}=\tilde{\lambda}_0.$$
This is a contradiction provided $\tilde{\lambda}_0=\lambda_0$. Here and the lemma below, we pre-assume that $\tilde{\lambda}_0=\lambda_0$ which will be proved soon after Lemma \ref{L1.1.4}. We will also mention in the proof that the pre-assumption is reasonable.
$\qquad\Box$

\begin{lem}\label{L1.1.3} The a.e. eigenfunction $g$  of $\lambda_{0}$ is either positive or negative everywhere.
\end{lem}

{\it Proof}\quad\rm If $\lambda_0=0$, then $g$ must be a constant and so the assertion is obvious. Now, let $\lambda_{0}>0$.
By Lemma \ref{L1.1.2}, without loss of generality, assume that $g'|_{(0,D)}<0$ and $g(0)>0$
\footnote{About $g(0)>0$. Since $Lg=-\lambda_{0}g$ on $(0,D)$, by (\ref{ff1}), we have ${\d g(s)}/{\d\nu}= -\lambda_0 \int_0^s g \d \mu, s\in (0,D)$.
So $\int_{0}^{s}g\d\mu>0$ by $g'<0$. This implies that $g(0)>0$. Otherwise,
one would get a contradiction with $\int_{0}^{s}g\d\mu>0$ since $g'<0$ and then $g\leqslant 0$.}
. We need only to prove that $g\ne 0$ on $(0, D)$. If otherwise $g(p)=0$ for some $p\in (0, D)$, then, since $\lambda_0^{(0,p)}$ is the minimal ND-eigenvalue on $(0,p)$, the eigenequation restricted to $(0, p)$ shows that
\footnote{The proof of $\lambda_{0}\geqslant\lambda_{0}^{(0,p)}$ is similar to the one of $\lambda_{0}\geqslant\lambda_1^{(0,p)}$ given in the last proof.}
$$\lambda_{0}\geqslant \lambda_0^{(0, p)}> \lambda_0^{(0, D)}=\tilde\lambda_0=\lambda_0,$$
which is a contradiction. $\qquad \Box$

Because of \eqref{ff1}, we have $I(g)^{-1}\equiv \lambda_0$. This explains where the operator $I$ comes from. Next, from \eqref{fnew2}, we
have
\begin{equation}\label{ff2}g(x)-g(D)=\lambda_0 \int_x^D \nu(\d s)\int_0^s g\d\mu.\end{equation}
When $D<\infty$, since $g(D)=0$ by our boundary condition, we
obtain  $I\!I(g)^{-1}\equiv \lambda_0$. This explains the meaning
of the operator $I\!I$. To show that the last assertion holds
even for $D=\infty$, it is necessary to prove that $g(\infty)=0$. This is impossible if $\lambda_0=0$ since then $g$ can be an arbitrary non-zero constant.

\begin{lem}\label{L1.1.4} Let $D=\infty$. If $\lambda_0>0$, then
its a.e. eigenfunction $g$ satisfies $g(\infty)=0$. \end{lem}

{\it Proof}\quad\rm Without loss of generality, by Lemmas \ref{L1.1.2} and \ref{L1.1.3},
assume that $g'|_{(0, D)}<0$ and $g|_{[0,D)}>0$.

(a) By what we have just seen and the decreasing property of $g$, we have
$$\frac{g(x)-g(\infty)}{\lambda_0}= \int_x^\infty \nu(\d s)\int_0^s g\d\mu\geqslant g(\infty) \int_x^\infty \nu(\d s)\int_0^s \d\mu.$$
Thus, $g(\infty)=0$ once $\int_0^\infty \nu(\d s)\int_0^s \d\mu=\infty$ (which is the uniqueness criterion for the semigroup or the nonexplosive criterion for the minimal process) since the left-hand side is finite.

(b) Otherwise, we have
$$M(x):=\int_x^\infty \nu(\d s)\int_0^s \d\mu<\infty,\qquad x\in (0, D).$$
Let $f=g-g(\infty)$ and suppose that $g(\infty)>0$. Then $f\in\mathscr{F}_{I\!I}$ and moreover,
$$\begin{aligned}
fI\!I(f)(x)={\lambda_{0}}^{-1}{\big(g(x)-g(\infty)\big)}-g(\infty)M(x)
=\lambda_0^{-1}f(x)-g(\infty)M(x). \end{aligned}$$
We arrive at
$$\sup_{x\in(0, \infty)}I\!I(f)(x)=\frac{1}{\lambda_{0}}-g(\infty)\inf_{x\in(0,\infty)}\frac{M(x)}{f(x)}. $$
Since $f(\infty)=0$ and $M(\infty)=0$, by Cauchy's mean value theorem, we have
$$\begin{aligned}
\inf_{x\in (0,\infty)}\frac{M(x)}{f(x)}
     &\geqslant\inf_{x\in (0,\infty)}\frac{M'(x)}{f'(x)}
    =\inf_{x\in (0,\infty)}-
    \frac{e^{-C(x)}}{g'(x)} \int_{0}^{x}\frac{e^{C(u)}}{a(u)}\d u\\
    &\geqslant\inf_{x\in (0,\infty)}-\frac{e^{-C(x)}}{g(0)g'(x)}\!\int_{0}^{x}g\,\d\mu  \quad(\text{since } g'\!<\!0\text{ and } g\!>\!0\text{ on }(0,D))
    \\&=\inf_{x\in (0,\infty)}\frac{1}{g(0)}I(g)(x)
    =\frac{1}{\lambda_{0}g(0)}>0.\end{aligned}$$
Inserting this into the previous equation, it follows that
$\lambda_0\!\!<\!\inf_{x\in(0,\infty)}\!\!I\!I(f)(x)^{-1}$.
But $\inf_{x\in(0, \infty)}I\!I(f)(x)^{-1}\leqslant \lambda_0 $
is a part of Theorem \ref{th1}\,(2) and will be proved soon below, without using the properties of the a.e. eigenfunction $g$. We have thus obtained a contradiction.  $\qquad \Box$

{\it Proof of Theorem $\ref{th1}$ and Proposition $\ref{p11.1.1}$}\quad
Similar to the proof of \rf{r2}{Theorem 2.4 and Proposition 2.5}, we can prove the assertions by two circle arguments.

 To prove the lower estimates, we adopt the following circle arguments:
\footnote{The details are given in Appendix A.1.}
\begin{eqnarray}
\lambda_{0}\hskip-0.6cm&&\geqslant\tilde{\lambda}_{0}\geqslant\sup\limits_{f\in\mathscr{F}_{I}}\inf\limits_{x\in(0,D)}I(f)(x)^{-1}\nonumber\\
&&\qquad\, =\sup\limits_{f\in\mathscr{F}_{I}}\inf_{x\in(0,D)}I\!I(f)(x)^{-1}
=\sup\limits_{f\in\mathscr{F}_{I\!I}}\inf\limits_{x\in(0, D)}I\!I(f)(x)^{-1}\label{fn4}\\
&&\geqslant\sup\limits_{h\in\mathscr{H}}\inf\limits_{x\in(0, D)}R(h)(x)\geqslant\lambda_{0}.\label{fn5}\end{eqnarray}
For the upper estimation of $\lambda_{0}$, we adopt the following circle arguments:
 \footnote{The details are given in Appendix A.2.}
\begin{eqnarray}
\lambda_{0}\hskip-0.6cm&&\leqslant\inf_{f\in\mathscr{{\widetilde F}}_{I\!I}\bigcup\mathscr{{\widetilde F}}_{I\!I}'}
\,\sup\limits_{x\in \supp(f)}I\!I(f)(x)^{-1}
=\inf\limits_{f\in\mathscr{{\widetilde F}}_{I\!I}}\,\sup\limits_{x\in
\supp(f)}I\!I(f)(x)^{-1}\label{fn6}\\
&&= \inf\limits_{f\in\mathscr{{\widetilde F}}_{I}}\,\sup\limits_{x\in
\supp(f)}I\!I(f)(x)^{-1}=
\inf\limits_{f\in\mathscr{{\widetilde F}}_{I}}\,\sup\limits_{x\in(0, D)}I(f)(x)^{-1}\nonumber\\
&&\qquad\qquad=\inf\limits_{f\in\mathscr{{\widetilde F}'}_{I}}\,\sup\limits_{x\in(0,
D)}I(f)(x)^{-1}\label{fn7}\\
&&\leqslant \inf\limits_{h\in\mathscr{\widetilde{H}}}\,\sup\limits_{x\in(0,
D)}R(h)(x)\leqslant\lambda_{0}\label{fn8}.
\end{eqnarray}

\rm  In fact,  most parts of the proof here are parallel to those in the discrete case (see \rf{r2}{Section 2}).
Actually, one can follow the cited proofs with some changes illustrated here. For instance, to prove
$\tilde{\lambda}_{0}\geqslant\sup_{f\in\mathscr{F}_{I\!I}}\inf_{x\in(0,
D)}I\!I(f)(x)^{-1}$,
following  \rf{r2}{Part I (a) of the proof of Theorem 2.4 and Proposition 2.5},
let $g$ (irrelated to the eigenfunction) be a test function of $\tilde{\lambda}_{0}$: $g\in {\scr C}^{1}(0, D)\cap {\scr C}[ 0, D]$, $g(D)=0$, and $\mu\big(g^{2}\big)=1$. Then for every $h$ with $h|_{(0,D)}>0$, we have
$$\begin{aligned}
1&=\mu(g^{2})=\int_{0}^{D}\frac{e^{C(x)}}{a(x)}\bigg(\int_{x}^{D}g'(t)\d t\bigg)^{2}\d x\\
 &\leqslant\!
\int_{0}^{D}\!\frac{e^{C(x)}}{a(x)}\d x\!\!\int_{x}^{D}\!\frac{e^{C(t)}}{h(t)}g'(t)^2\d t\!\!\!\int_{x}^{D}\!\frac{h(s)}{e^{C(s)}}\d s
\;\;\text{(by Cauchy-Schwarz's inequality)}\\
&=\int_{0}^{D}\frac{e^{C(t)}}{h(t)}g'(t)^2\d t\int_{0}^{t}\frac{e^{C(x)}}{a(x)}\d x\int_{x}^{D}\frac{h(s)}{e^{C(s)}} \d s
\;\;\text{(by Fubini's Theorem)}\\
&\leqslant
 D(g)\sup\limits_{t\in(0, D)}\frac{1}{h(t)}\int_{0}^{t}\frac{e^{C(x)}}{a(x)}\d x\int_{x}^{D}\frac{h(s)}{e^{C(s)}}\d s\\&
=:D(g)\sup\limits_{t\in(0, D)}H(t).\end{aligned}$$
For $f\in
\mathscr{F}_{I\!I}$ satisfying $\sup_{x\in(0, D)}I\!I(f)(x)<\infty$, we specify
$h(t)=\int_{0}^{t}{a(s)}^{-1}{e^{C(s)}}\\f(s) \d s$. Then
by Cauchy's mean value theorem, it follows that
$$\begin{aligned}
\sup\limits_{t\in(0, D)}H(t)\leqslant\sup\limits_{x\in(0, D)}\frac{1}{f(x)}\int_{x}^{D}e^{-C(s)}\d s\int_{0}^{s}\frac{e^{C(u)}}{a(u)}f(u)\d u=\sup\limits_{x\in(0, D)}I\!I(f)(x).\end{aligned} $$
Hence,
$$\inf_{x\in(0, D)}I\!I(f)(x)^{-1}\leqslant\inf_{t\in(0, D)}H(t)^{-1}\leqslant D(g).$$
Making infimum with respect to $g$, we obtain the required assertion. We have also completed the proofs of Lemmas \ref{L1.1.2}--\ref{L1.1.4}.

From now on in this section, we assume that the a.e. eigenfunction (say $g$)  satisfies $g>0$ and $g'<0$ on $(0, D)$, $g'(0)=0$, and $g(D)=0$ (recall that $g(D)=\lim_{x\to D} g(x)$ if $D=\infty$).

As mentioned before Lemma \ref{L1.1.4}, the operators $I$ and $I\!I$ are all come from the eigenequation. Here we show that
so is the operator $R$. Rewrite the eigenequation as
$$-\frac {Lg}{g}=\lambda_0$$
which is meaningful since $g>0$. To simplify the left-hand side,
in the discrete case, one uses the ratio $g(x+1)/g(x)$. However,
this is useless in the present continuous situation. What instead is using the function $h=g'/g$. Then
$$-\frac {Lg}{g}=-(a h^2+ bh +a h')=R(h).$$
The conditions $g>0$ and $g'<0$  on $(0,D)$ lead to the restraint
$h|_{(0,D)}<0$ in defining ${\scr H}$.
Note that the inverse transform $h\to g$ is unique up to a
positive constant:
$$g(x)= \exp\bigg[\int_0^x h(u) \d u\bigg].$$
The restraint allowing $h=0$ in the definition of ${\scr H}$
is to include the degenerated case that $g'\equiv 0$ when $\lambda_0=0$ (then $D=\infty$ by hypothesis \eqref{f1.0}).
Clearly, the use of $R$ is essentially the use of $L$. Since
this, we make the continuous condition on $a$ and $b$
once concerning with $R$.
Because of this point, we need two additions in the above circle
arguments: the right-hand side of \eqref{fn4} is not less than $\lambda_0$ and
the right-hand side of \eqref{fn7} is no more than $\lambda_0$. This is rather easy since for the a.e eigenfunction $g$, we have $I(g)^{-1}\equiv \lambda_0$ and $I\!I(g)^{-1}\equiv \lambda_0$ by \eqref{ff1}, \eqref{ff2} and Lemma \ref{L1.1.4}. Actually, the required assertion was also contained in the corresponding proof of  the discrete
situation.

As another illustration of the proof moving from the discrete case to the continuous one, we consider a proof for the upper estimates. For instance, we prove that
$$\lambda_{0}\leqslant\inf_{f\in\mathscr{\widetilde{F}}_{I\!I}\bigcup\mathscr{\widetilde{F}}_{I\!I}'}
\sup_{x\in \supp\,(f)}I\!I(f)(x)^{-1}.$$

Before  moving to the details, let us mention that,  for the upper estimates of $\lambda_{0}$, we are actually using a comparison between $\lambda_{0}$ and $\lambda_{0}^{(0, x_{0})}$. Thus, for the upper estimates of $\lambda_{0}$, we indeed use the restriction on $[0, x_{0}]$ of the test functions, ignoring the behavior of them out of $[0, x_{0}]$.

Given $f\in\widetilde{\mathscr{F}}_{I\!I}$ with $f=f\mathbbold{1}_{[0, x_{0})}$ for some $x_0\in (0, D)$,  let $$g=fI\!I(f)
\mathbbold{1}_{\supp\,(f)}.$$  Then $g\in L^{2}(\mu)$.
Since
$$\Big[e^C g'\Big](x)=-\int_{0}^{x}f\,\d\mu \qquad\text{on $[0,x_{0})$},$$
by the integration by parts formula, we have
$$\aligned D(g)&=\int_{0}^{D}e^{C(x)}g'(x)^{2}\d x
=-\int_{0}^{x_0 } \big(\d g(x)\big)\int_{0}^{x}f\,\d\mu\\
&=-\int_0^{x_{0}}f(t)\int_t^{x_{0}}\d g(t)\mu(\d t)\quad(\text{by Fubini's Theorem})\\
&=\int_{0}^{x_0 } f(t) (g(x_{0})-g(t))\mu(\d t)\\
&=\int_{0}^{x_0 } f g \d\mu\quad (\text{since $g(x_{0})=0$}).\endaligned$$
Hence
$$\aligned
   D(g)&
   \leqslant\int_{0}^{x_{0} }g^{2}\,\d\mu\sup\limits_{(0, x_{0} )}\frac{f}{g}=\mu(g^{2})\sup\limits_{x\in(0, x_{0})}I\!I(f)(x)^{-1}.
   \endaligned $$
Since $g\in L^{2}(\mu)$, it follows that
\begin{equation}\label{f1.2}\lambda_{0}\leqslant\frac{D(g)}{\mu(g^2)}\leqslant\sup_{\rm supp\,(f)}I\!I(f)^{-1}.\end{equation}
for  every $ f\in\widetilde{\mathscr{F}}_{I\!I}$.
It remains to show that the same assertion holds for every $f\in \mathscr{\widetilde{F}}_{I\!I}'$.
Recall that in the proof above, the conclusion $g\in L^{2}(\mu)$
comes from the finiteness of $x_0$. Otherwise, if $x_{0}=D=\infty$, then $ f\in\widetilde{\mathscr{F}}_{I\!I}'$ means that the function $g=fI\!I(f)$ is assumed to be in
$L^{2}(\mu)$, and so the proof above still works. So we obtain again the required assertion.

Hopefully, we have explained enough the difference between the discrete and the continuous cases. Now, one may follow \rf{r2}{Proof of Theorem 2.4 and Proposition 2.5} (quite long and technical) to complete the whole proof. $\qquad \Box$

Before moving further, let us mention a fact about the localizing procedures used in Theorem \ref{t11.2.2}\,(2).
Instead of the approximating to the infinite state space
($D=\infty$) by finite ones, it seems more natural to use the truncating procedure for the test function $f$: $f^{(n)}=f \mathbbold{1}_{[0, x_{n})}$ with $x_n\uparrow \infty$. The next result shows that such a procedure is not practical in general.

\begin{rem}\label{re1.1.3} Assume that hypothesis \eqref{f1.0} holds. Let $D=\infty$ and $g$ be the eigenfunction of
$\lambda_{0}>0$,\ define $g^{(n)}=g \mathbbold{1}_{[0, x_{n})}$ for some $x_{n}\in(0, \infty)$. Then
\begin{equation*}\inf\limits_{x\,\in\,
\supp\,(g^{(n)})}I\!I(g^{(n)})(x)=0. \end{equation*} In particular,\
$\inf_{x\,\in\, \supp\,(g^{(n)})}I\!I(g^{(n)})(x)$ does not converge to
$\lambda_{0}$ as $x_{n}\to \infty$. \end{rem}
{\it Proof}\quad\rm By the definition of $g^{(n)}$, we have
$$\aligned
   \inf\limits_{x\,\in\, \supp\,(g^{(n)})}I\!I(g^{(n)})(x)&=\inf_{x\,\in\,[0, x_{n})}\frac{1}{g^{(n)}(x)}\int_{x}^{x_{n}}e^{-C(s)}\d s\int_{0}^{s}g^{(n)}\,\d\mu\\
   &=\inf_{x\,\in\,[0, x_{n})}\frac{1}{g(x)}\int_{x}^{x_{n}}e^{-C(s)}\d s\int_{0}^{s}g\,\d\mu\\
   &=\inf_{x\,\in\,[0, x_{n})}\frac{1}{g(x)}\int_{x}^{x_{n}}\big(-\lambda_{0}^{-1}g'(s)\big)\d s\qquad\big(\text{by }\eqref{ff1}\big)\\
   &=\inf_{x\,\in\,[0, x_{n})}\frac{1}{\lambda_{0}g(x)}\big(g(x)-g(x_{n})\big)
 \\&=\inf_{x\,\in\,[0,x_{n})}\frac{1}{\lambda_{0}}\bigg(1-\frac{g(x_{n})}{g(x)}\bigg)\\&=0
 \qquad(\text{since } g\in {\scr C}[0,D]\text{ and } g \downarrow\downarrow).
  \quad\Box\endaligned
  $$

{\it Proof of Proposition $\ref{t11.1.3}$}\quad\rm $(1)$ Let $g\in {\scr C}^{1}(0,D)$ with $g>0$ and $g'<0$ on $(0, D)$, and  let $\bar{h}(x)=-e^{-C(x)}\int_{0}^{x}g \,\d\mu$. Then
\footnote{$\bar{h}(0)=0,\,\bar{h}'=-bh/a-g/a<-bh/a.$ So $\bar{h}\in{\scr C}^{2}(0,D)\cap{\scr C}[0,D]$. Thus, $\bar{h}\in\mathscr{H_{\ast}}$. Moreover, $a\bar{h}''+b\bar{h}'=-g<0$.}
 $\bar{h}\in\mathscr{H_{\ast}}$ and
\begin{equation*}
\overline R(\bar{h})(x)=-\frac{(a\bar{h}'+b\bar{h})'(x)}{\bar{h}(x)}=\frac{g'(x)}{\bar{h}(x)}>0.
\end{equation*}
This clearly implies that $\sup_{h\in\mathscr{H}_{\ast}}\inf_{x\in(0,D)}\overline
R(h)(x)\geqslant 0$.

$(2)$ Without loss of generality, assume that $\lambda_0>0$. Since $a,\, b\in {\scr C}^1(0,D)$, there exists an eigenfunction $g$ such that
\footnote{We have $g'(0)=0$, $g\in{\scr C}^{2}[0, D]$, $g>0$ and $g'<0$ on $(0,D)$ by Lemmas \ref{L1.1.0}, \ref{L1.1.2} and \ref{L1.1.3}. Since $\bar{h}=g'$ and $a, b\in{\scr C}[0,D]\cap{\scr C}^1(0,D)$, we can see that $\bar{h}(0)=0$  and
\begin{equation*}h'=g''=-\frac{\lambda_0 g+b g'}{a}\in{\scr C}^1(0,D),\quad a\bar{h}'+b\bar{h}=ag''+bg'=-\lambda_{0}g<0. \end{equation*}
So $\bar{h}\in{\scr C}^2{(0,D)}\cap{\scr C}[0,D]$,\, $\bar{h}'<-{a}^{-1}{b}\bar{h}$ and then $\bar{h}\in \mathscr{H_{\ast}}$.}
 $\bar{h}:=g'\in\mathscr{H}_{\ast}$ and
$$\aligned
\overline R(\bar{h})(x)=-(Lg)'(x)/g'(x)\equiv\lambda_{0}.
\endaligned$$
Thus
 $$\aligned\sup_{h\in\overline{\mathscr{H}}}\inf_{x\in(0,D)}\overline R(h)(x)\geqslant\sup_{h\in\mathscr{H}_{\ast}}\inf_{x\in(0,D)}\overline R(h)(x)\geqslant\lambda_{0}.\endaligned$$

 Now, one can complete the proof following that in the discrete case (\rf{r2}{Proof of Proposition 2.7})
 \footnote{The details are given in Appendix A.3.}
 .$\qquad \Box$

To prove Theorem \ref{t11.2.1}, we need the following result.

\begin{lem}\label{l00} Given two nonnegative, measurable, and
locally integrable functions $m$ and $n$ on $[0, D]$,  suppose that
$$\int_{0}^{D}n(y)\d y<\infty\quad \text{\rm  and}\quad c:=\sup\limits_{x\in(0,D)} \int_{0}^{x}m(y)\d y\int_{x}^{D}n(y)\d y<\infty.$$
Set $\psi(x)=\int_{x}^{D}n(y)\d y$.
Then for every $r\in (0, 1)$, we have
\begin{equation*}\int_{0}^{x}m(y)\psi^r(y)\d y\leqslant\frac{c}{1-r}\psi^{r-1}(x),\qquad  x\in(0, D).
\end{equation*}
\end{lem}

{\it Proof}\quad
\footnote{The details are given in Appendix A.4.}
\rm Let $M(x)=\int_{0}^{x}m(y)\d y$. Noticing that  $M'(x)=m(x)$ and $M\psi\leqslant c$,
we obtain the assertion  by using the integration by parts formula. $\qquad \Box$

{ \it Proof of Theorem $\ref{t11.2.1}$}\quad\rm To prove the lower estimate, without loss of generality, assume that $\delta<\infty$.
Applying Lemma \ref{l00} to $m(x)={e^{C(x)}}/{a(x)}$ and $n(x)=e^{-C(x)},$ we get
\begin{equation*}\int_{0}^{x}\varphi^{r}(y)\mu(\d y)=\int_{0}^{x}\varphi^{r}(y)m(y)\d y\leqslant\frac{\delta}{1-r}\varphi^{r-1}(x),\qquad x\in(0, D). \end{equation*}
Put $f=\varphi^{r}$. Then $f\in\mathscr{F}_{I}$ and
$I(f)(x)\leqslant{\delta}/{(r-r^{2})}.$ Optimizing the inequality with respect to $r$, it follows that
\begin{equation} I(f)(x)\leqslant\inf_{0<r<1}{\delta}/{(r-r^{2})}=4\delta \label{f1.1.16}.\end{equation}
We have thus proved the lower estimate.

For the upper estimate
\footnote{The details are given in Appendix A.5.}
 , we choose the test function
$f=\nu (x_{0}\vee\cdot, x_{1}) \mathbbold{1}_{[0, x_{1})}$ for some $ x_{0}$, $x_{1}\in[0, D)$ with $x_0<x_1$. Then, the assertion follows by using either the variational formula
for upper estimate given by Theorem \ref{th1}\,(1)
$$\lambda_{0}\leqslant\inf_{f\in\mathscr{{\widetilde F}}_{I}}\sup_{x\in(0, D)}I(f)(x)^{-1}$$
or the classical variational formula:
 $$\lambda_{0}^{-1}=\tilde{\lambda}_{0}^{-1}=\lambda_{\ast}^{-1}\geqslant\!\!\sup_{x_{0},x_{1}:x_{0}<x_{1}}
\frac{\big\|f_{1}^{x_{0},x_{1}}\big\|}{D\big(f_{1}^{x_{0},x_{1}}\big)}$$
and then letting $x_{1}\rightarrow D$.

At last, if $\nu(0,D)=\infty$, then we have $\nu(x,D)=\infty$ because of hypothesis \eqref{f1.0}. Furthermore, $\mu(0,x)\nu(x,D)=\infty$ for every $x\in(0,D)$. So $\delta=\infty$ and $\lambda_{0}=0$. If $\int_0^\infty\mu(0, x)\,\nu(\d x)<\infty$, then for each $x\in(0,D)$, we have $$\aligned\mu(0,x)\nu(x,D)=\int_{x}^{D}\mu(0,x)\nu(\d t)<\int_{x}^{\infty}\mu(0,t)\nu(\d t)<\int_0^\infty\mu(0, x)\,\nu(\d x)<\infty\endaligned.$$ Hence, $\delta<\infty$ and $\lambda_{0}>0$.$\qquad \Box$

{\it Proof of Theorem $\ref{t11.2.2}$ and Corollary $\ref{cor1}$}\quad
\footnote{The details are given in Appendix A.6.}
\rm Simply follow \rf{r2}{Proof of Theorem 3.2 and Corollary 3.3}. We mention that the proof of ``$\delta_{1}'\leqslant2\delta$'' and the computation of
$\delta_{1}'$ are not easy.$\qquad \Box$
\noindent \\[4mm]

\noindent{\bbb 4\quad The DN-case }\\[0.1cm]
We now turn to study the DN-case.
As in Section 2, we use the same notation ${\scr C}[0, D],\ {\scr C}^k(0, D)$ and the operator $L$.
 The main different point for the eigenequation $L g=-\lambda_{0} g$ is the boundary conditions: $g(0)=0$ and $g'(D)=0$ if $D<\infty$. Now define
\begin{equation}\lambda_{0}\!=\inf\bigg\{\frac {D(f)}{\mu\big(f^{2}\big)}\!: f\!\in\! {\scr C}^1(0,D)\cap {\scr C}[0,D],\,D(f)\!<\!\infty, f(0)\!=\!0,\, f\!\ne\! 0\bigg\},\label{1f1.1}\end{equation}
where
$$\aligned &D(f)=\int_{0}^{D}a {f'}^2\d\mu,\quad &\mu(\d x)=\frac{e^{C(x)}}{a(x)}\d x,\quad C(x)=\int_{0}^{x}\frac{b(u)}{a(u)} \d u.
\endaligned$$
Again, define $\nu(\d x)=e^{-C(x)}\d x$.
Here, we have used the hypothesis \eqref{f1.0}.
The restraint ``$D(f)<\infty$'' in \eqref{1f1.1} is to avoid
$\infty/\infty$ since we allow $\mu\big(f^{2}\big)=\infty$.
Then the restraint ``$f\ne 0$'' is needed to avoid $0/0$.
Note that the restriction on the set ${\scr C}_K$ of test functions disappears in \eqref{1f1.1}. This means that the maximal Dirichlet form or the maximal process is used here, instead of the minimal one used in Section 2. In other words, we do not assume the uniqueness of the semigroup, which is different from what we studied earlier in [\ref{r4}\,--\,\ref{r6}] and \cite{r1}. The constant $\lambda_{0}$ defined above describes the optimal constant $C=\lambda_{0}^{-1}$ in the following $weighted\ Hardy \ inequality:$
$$\aligned
\mu\big(f^{2}\big)\leqslant C D(f),\quad f(0)=0.
\endaligned$$
(See \rf{r6}{Section 5.2}). In other words, we are studying the weighted Hardy inequality in this section. To save the notation, we use the same notation $\lambda_{0},$ $I,\ I\!I,\ R $ and so on as before, each of them plays a similar role but may have different meaning in different context.

Before going to our main text, we note that in the definition of $\lambda_{0}$, one may replace $\mathscr{C}^{1}(0,D)\cap{\scr C}[0,D]$ by $\mathscr{A}[0,D]$ as shown by Lemma \ref{new}.

Now, we review some notation defined originally in \cite{r5, r1}  and introduce some new ones as follows.
$$\aligned
&I(f)(x)=\frac{e^{-C(x)}}{f'(x)}\int_{x}^{D}f \,\d\mu \qquad \mbox{(single integral form)},\\
&I\!I(f)(x)=\frac{1}{f(x)}\int_{0}^{x}\nu(\d s)\int_{s}^{D}f \, \d\mu\qquad \mbox{(double integral form)},\\
&R(h)(x)=-(ah^{2}+bh+ah')(x)\qquad \mbox{(differential form)}.
\endaligned$$
The domains of $I$, $I\!I$ and $R$, respectively, are as follows.
$$\aligned
&\mathscr{F}_{I}=\{f: f\in {\scr C}^1(0, D)\cap {\scr C}[0, D], f(0)=0, \text{ and } f'|_{(0,D)}>0\},\\
&\mathscr{F}_{I\!I} =\{f: f\in {\scr C}[0,D], f(0)=0, \text{ and } f|_{(0, D)}>0\},\\
&\mathscr{H}=\Big\{h: h\in {\scr C}^1(0, D)\cap {\scr C}[0, D],\,h|_{(0, D)}>0,\text{ and }\mbox{$\int_{0+}\!\!h(u)\d u=\infty$}\Big\},
\endaligned$$
where $\mbox{$\int_{0+}$}$ means $\int_0^{\varepsilon}$ for sufficiently small $\varepsilon>0$. These sets
 \footnote{About the set of $\mathscr{H}$: The functions $g$ and $h$ are one to one: $h=g'/g$ and conversely
 $g(x)=g(\varepsilon)\exp\big[\int_{\varepsilon}^{x}h(u)\d u\big],\;x\in(0,D)$. So $g(0)=0$ implies $h(0)=\infty$
 and conversely $g(0)=0$ is implied by \mbox{$\int_{0+}h(u)\d u=\infty$}.}
 are used for the estimates on lower bounds of $\lambda_0$. For the upper bounds, we have the following domains.
$$\aligned&\mathscr{{\widetilde F}}_{I}=\big\{f: \exists x_{0}\in(0, D),\,f\in {\scr C}^1(0, x_{0})\cap {\scr C}[0, D],\,f(0)=0, f=f(\cdot\wedge x_{0}),\\
&\text{\hskip4em\; and } f'|_{(0,\, x_{0})}\!>0\big\},\\&\mathscr{{\widetilde F}}_{I\!I}=\big\{f: \exists x_{0}\in(0, D),\,f\in {\scr C}[0, x_{0}],\,f(0)=0,\,f=f (\cdot\wedge x_{0})\text{ and }f|_{(0, x_0)}>0 \big\},\\
&\mathscr{\widetilde{H}}=\Big\{h:  \exists x_{0}\in(0, D),\,h\in {\scr C}^1(0,x_{0})\cap {\scr C}[0,D],\,h|_{(0, x_0)}>0,\;\mbox{$\int_{0+}h(u)\d u=\infty$},\\
&\text{\hskip4em\;} \,h|_{[x_0, D]}=0,\text{ and } \sup_{(0,  x_0)}\big(ah^{2}+bh+ah'\big)<0\Big\}.
\endaligned$$
Besides, we need also
$$\aligned
\mathscr{{\widetilde F}'}_{I\!I}=\big\{f: f>0,\,f\in {\scr C}[0,D], \text{ and } fI\!I(f)\in L^{2}(\mu)\big\}.
\endaligned$$

Under  hypothesis \eqref{f1.0}, if $\mu(0,D)=\infty$, then $\lambda_{0}$ defined by \eqref{1f1.1} is trivial.
Indeed, let $$f=\mathbbold{1}_{(\delta,\, D]}+h\mathbbold{1}_{[0, \delta]},$$
 where $h$ is chosen so that $h(0)=0$ and $f\in {\scr C}^1(0,D)\cap {\scr C}[0,D]$ (For example, $h(x)=-x^{2}/\delta^{2}+2x/\delta$).
Then $D(f)\in (0,\infty)$ and $\mu\big(f^{2}\big)=\infty$.
It follows that  $\lambda_{0}=0.$

Otherwise, $\mu(0,D)<\infty$. Then for every $f$ with $\mu\big(f^{2}\big)=\infty,$ by setting $f^{(x_{0})}=f(\cdot\wedge x_{0})\in L^{2}(\mu)$,  we have
$$
\aligned\infty>D\big(f^{(x_{0})}\big)\uparrow D(f)\quad \text{and}\quad  \infty>\mu\Big({f^{(x_{0})}}^{2}\Big)\rightarrow
\mu\big(f^{2}\big)\qquad\text{as}\; x_{0}\rightarrow D.
\endaligned$$ In other words, for each non-square-integrable function $f,$ both $\mu\big(f^{2}\big)$ and $D(f)$ can be approximated by a sequence of square-integrable ones. Hence,
we can rewrite $\lambda_{0}$ as follows.
\begin{equation}
\lambda_{0}=\inf\big\{D(f):  \mu\big(f^{2}\big)=1, f(0)=0, \text{ and } f\in {\scr C}^1(0,D)\cap {\scr C}[0,D]\big\}.\label{1f1.2}
\end{equation}
In this case, as will be seen soon but not obvious, we also have  $$\aligned
\lambda_{0}=\inf\big\{D(f):\, &\mu\big(f^{2}\big)=1,\, f(0)=0,\,f=f(\cdot\wedge x_{0}), \,f\in {\scr C}^1(0,x_{0})\cap {\scr C}[0,x_{0}]\\&\text{for } \text{some } x_{0}\in (0,D)\big\}=:\tilde{\lambda}_{0}.\endaligned$$

Now we introduce our main results. Their relations are
very much the same as indicated in Section 2, except that the test function used in Theorem \ref{1t1.2} is $\nu(0, x)^{\gamma}$
but not $\nu(x, D)^{\gamma}$\,$(\gamma=1/2$ or $1)$.

\begin{thm}\label{1t1.1}
\footnote{The details of the proof are given in Appendix B.1.}
Let hypothesis \eqref{f1.0} hold. Assume that $\mu(0,D)<\infty$. Then $\lambda_0$ defined by \eqref{1f1.1} or \eqref{1f1.2} coincides with $\tilde{\lambda}_{0}$ and
the following variational formulas hold.

$(1)$ Single integral form:
$$
\aligned
\inf\limits_{f\in\mathscr{{\widetilde F}}_{I}}\sup\limits_{x\in(0, D)}I(f)(x)^{-1}=\lambda_{0}=\sup\limits_{f\in\mathscr{F}_{I}}\inf\limits_{x\in(0,D)}I(f)(x)^{-1}.
\endaligned$$

$(2)$ Double integral form:
$$\aligned
\lambda_{0}&=\inf_{f\in\mathscr{{\widetilde F}}_{I}}\sup_{x\in
    (0,D)}\!\!I\!I(f)(x)^{-1}
   \!\!=\!\!\inf_{f\in\mathscr{{\widetilde F}}_{I\!I}}\sup_{x\in
    (0,D)}\!I\!I(f)(x)^{-1}
    \!=\!\!\inf_{f\in\mathscr{{\widetilde F}}_{I\!I}\cup \mathscr{{\widetilde F}}_{I\!I}'}\sup_{x\in
    (0,D)}\!\!I\!I(f)(x)^{-1}\\ \lambda_0&=\sup\limits_{f\in\mathscr{F}_{I}}\inf\limits_{x\in(0, D)}I\!I(f)(x)^{-1}
=\sup\limits_{f\in\mathscr{F}_{I\!I}}\inf\limits_{x\in(0, D)}I\!I(f)(x)^{-1}.
\endaligned$$
Moreover, if $a,\ b\in {\scr C}[0,D],$ then we also have a

$(3)$ differential form:
$$\aligned
 \inf\limits_{h\in\mathscr{\widetilde{H}}}\sup\limits_{x\in(0,
D)}R(h)(x)=\lambda_{0}=\sup\limits_{h\in\mathscr{H}}\inf\limits_{x\in(0, D)}R(h)(x).
\endaligned
$$ \end{thm}

\begin{thm}\label{1t1.2}
$(\text{\rm Criterion and basic estimates})$
Let hypothesis \eqref{f1.0} hold. Then $\lambda_{0}$ defined by \eqref{1f1.1} $($or equivalently $\tilde{\lambda}_{0}$ provided $\mu(0,D)<\infty)$ is positive iff $$\delta:=\sup\limits_{x\in(0, D)}\nu(0, x)\,\mu(x, D)<\infty. $$
More precisely, we have
$$(4\delta)^{-1}\leqslant\lambda_{0}\leqslant\delta^{-1}.$$
 In particular, we have $\lambda_{0}=0$ if $\mu(0,D)=\infty$, and $\lambda_{0}>0$ if
 $$D<\infty\quad\text{ or }\quad \int_{0}^{D}\Big({a(u)}^{-1}{e^{C(u)}}+e^{-C(u)}\Big)\,\d u<\infty.$$
\end{thm}
{\it Proof}\quad
\rm The result  was proved in \rf{r3}{Theorem 1.1} except the case that $\mu(0,D)=\infty$ in which case
$\lambda_{0}=0$ ($\delta=\infty$) and so the assertion is trivial.

\begin{thm}\label{1t1.3}
\footnote{The details of the proof  are given in Appendix B.2.}
$(\text{\rm Approximating procedure})$
Let hypothesis \eqref{f1.0} hold. Assume that $\mu(0,D)<\infty$ and $\delta<\infty$. Set $\varphi(x)=\nu (0, x)$ for $x\in (0, D)$.

$(1)$ Define $f_{1}=\sqrt{\varphi}$, $f_{n}=f_{n-1}I\!I(f_{n-1})$, $n\geqslant 2$, and  let $\delta_{n}=\sup\limits_{x\in(0, D)}I\!I(f_{n})(x)$, $n\geqslant 1$. Then $\delta_{n}$ is
decreasing in $n$ and
$$ \lambda_{0}\geqslant\delta_{n}^{-1}\geqslant({4\delta})^{-1},\qquad n\geqslant1.$$

$(2)$ For fixed $x_{0}\in (0, D)$, define
$$f_{1}^{(x_{0})}=\varphi(\cdot \wedge x_{0}),
\qquad f_{n}^{(x_{0})}=\big(f_{n-1}^{(x_{0})}I\!I\big(f_{n-1}^{(x_{0})}\big)\big)(\cdot \wedge x_0),\qquad n\geqslant2,$$
and let $\delta_{n}'=\sup_{x_{0}\in(0,D)}\inf_{x\in (0,D)}I\!I({f_{n}}^{(x_{0})})(x)$.
Then $\delta_{n}'$ is increasing in $n$ and
$$\delta^{-1}\geqslant
{\delta_{n}'}^{-1}\geqslant\lambda_{0},\qquad n\geqslant1.$$
Next, define
$$\bar{\delta}_{n}=\sup\limits_{x_{0}\in (0,D)}\frac{\big\|{f_{n}}^{(x_{0})}\big\|}{D\big({f_{n}}^{(x_{0})}\big)},\qquad n\geqslant 1. $$
Then ${\bar{\delta}_{n}}^{-1}\geqslant \lambda_{0}$, \, $\bar{\delta}_{n+1}\geqslant\delta_{n}'$\;  for every $n\geqslant 1$ and $\bar{\delta}_{1}=\delta_{1}$.
\end{thm}
\begin{cor}\label{1c1.1}
\footnote{The details of the proof  are given in Appendix B.3.}
 $(\text{\rm Improved estimates})$
We have the following estimates:
$$\aligned
\delta^{-1}\geqslant{\delta_{1}'}^{-1}\geqslant\lambda_{0}\geqslant\delta_{1}^{-1}\geqslant(4\delta)^{-1},
\endaligned$$
where
$$
\aligned
\delta_{1}&=\sup\limits_{x\in(0,D)}\frac{1}{\sqrt{\varphi(x)}}
\int_{0}^{D}\varphi(x\wedge \cdot)\sqrt{\varphi}\,\d\mu\\
&=\sup\limits_{x\in(0,D)}\bigg(\frac{1}{\sqrt{\varphi(x)}}
\int_{0}^{x}\varphi^{{3}/{2}}\d\mu
+\sqrt{\varphi(x)}\int_{x}^{D}\sqrt{\varphi}\,\d\mu\bigg),
\\
\delta_{1}'&=\sup\limits_{x\in(0,D)}\frac{1}{\varphi(x)}
\int_{0}^{D}\varphi(\cdot\wedge x)^{2}\d\mu\;\;\in[\delta,2\delta].
\endaligned$$
\end{cor}
\bigskip

Since the proofs of the above results are either known from
\cite{r3, r5} or parallel to \cite{r2}, here we make some remarks only.

\begin{rem}{\rm (1) As mentioned in \cite{r2}, the original proofs given in
\cite{r3, r5} are still suitable to support the idea using the maximal Dirichlet form instead of the uniqueness assumption.

(2) As discussed in the last section, it is natural to extend $a$ and $b$ from continuous to measurable in the case
using operators $I$ and $I\!I$ only.

(3) About the duality. Recall that
$$L=\frac{\d}{\d \mu} \frac{\d}{\d \nu}.$$
The dual operator of $L$ is simply defined as
$$L^*=\frac{\d}{\d \mu^*} \frac{\d}{\d \nu^*},
\qquad \mu^*:=\nu,\; \nu^*:=\mu.$$
For the boundaries, simply exchange the names Dirichlet and Neumann.
The basic results for these operators are
$\lambda_0(L)=\lambda_0(L^*)$ and $\delta=\delta^*$,
where $\lambda_0(L)$ and $\delta$ are defined in Section 2,
and $\lambda_0(L^*)$ and $\delta^*$ are defined in this section
replacing $L$ with $L^*$. The proof goes as follows.

(a) Reduce to finite $D$. By an approximating procedure we have
used many times before, it suffices to prove the assertion for
finite $D$. The point is that for $\lambda_0(L)$, one needs to
consider only the test functions having compact support; for $\lambda_0(L^*)$, it suffices to consider the test function
$f=f(\cdot \wedge x_0)$, where $x_0$ varies over $(0, D)$.

(b) By a standard smoothing procedure, one may assume that
$a$ and $b$ are smooth.

(c) The identity of $\lambda_0(L)$ and $\lambda_0(L^*)$ is a
combination of Proposition \ref{t11.1.3}\,(2) and Theorem \ref{1t1.1}\,(3). The discrete case was given in \rf{r2}{Section 5}. An alternative proof of this assertion is presented in
\cite{rn1} based on isospectral. Note that in the last proof,
the finiteness of $D$ is crucial, otherwise, the domains of
$L$ and $L^*$ are essential different unless the Dirichlet form corresponding to $L^*$ is assumed to be regular.

(4) When $D<\infty$, one may simply reverse the variable to obtain one from the other of the ND- and DN- cases. In this sense, the identity $\lambda_{0}(L)=\lambda_{0}(L^{\ast})$ stated in (3) is quite natural even though the duality is not a ``reverse transform''. When $D=\infty$, these two cases are certainly different since the Dirichlet boundary at 0 is touchable but not the one at $\infty$. We mention that the variational formulas and then the approximating procedure
in this section are different from those deduced by the
dual approach. It is interesting that in the discrete situation, the approximating procedure given by Theorem \ref{1t1.3} is often less powerful than those given by Theorem \ref{t11.2.2} in terms of duality. Similar phenomenon happens in the continuous situation as shown in \cite{rn2} with $D<\infty$.
}
\end{rem}
\noindent \\[4mm]

\noindent{\bbb 5\quad Supplement to the NN-case }\\[0.1cm]
Everything is the same as those in the last section except the mixed
eigenvalue $\lambda_0$ is replaced by
\begin{equation}
\lambda_1=\inf\big\{D(f): \mu(f)=0,\, \mu\big(f^{2}\big)=1,\; \text{ and } f\in {\scr C}^1(0,D)\cap {\scr C}[0,D]\big\}.
\end{equation}
Let us repeat that throughout this section, we
assume that hypothesis \eqref{f1.0} holds and $\mu(0,D)<\infty$.

The supplement consists of three parts. The first one is using
the maximal Dirichlet form instead of the uniqueness assumption of the semigroup. The second one is using the ``a.e. eigenfunction'' instead of ``eigenfunction''. These two parts have already been studied in the last section. See also \cite{r1} for some supplement to the
original paper. The third part is about the monotonicity of an
approximating procedure which we are going to study below.

Define
\begin{equation*}\bar{f}=f-\pi{(f)},\quad f_{1}=\sqrt{\varphi},\quad  f_{n}=\bar{f}_{n-1}I\!I\big(\bar{f}_{n-1}\big),\quad \eta_{n}=\sup\limits_{x\in(0, D)}I\big(\bar{f}_{n}\big)(x), \end{equation*}
where $\pi=\mu/\mu(0, D)$.
Here our main question  is about the monotonicity of $\{\eta_{n}\}$.
Unlike the sequences $\{\delta_n\}$ and $\{\delta_n'\}$ defined
in Theorems \ref{t11.2.2} and \ref{1t1.3}, their monotonicity is
simply twice applications of Cauchy's mean value theorem, the method does not work for the sequence $\{\eta_{n}\}$ since each
$\bar{f}_{n}$ can be zero in $(0, D)$.
We were unable to solve this problem for years until the appearance of the recent paper \rf{r2}{Section 6}, in which the problem was solved in the discrete context. Note that
$\lambda_1>0$ iff
$$\delta:=\sup\limits_{x\in(0, D)}\nu(0, x)\,\mu(x, D)<\infty$$ by \rf{r3}{Theorem 3.7}, \rf{r2}{Theorem 6.2}, and Theorem \ref{1t1.2}.

\begin{pro}\label{t11.2.3} Let hypothesis \eqref{f1.0} hold and assume that
$\delta <\infty$. Then the sequence $\{\eta_{n}\}$ defined above $\big($i.e. $\{\eta_n''\}$ in {\rm\rf{r5}{Theorem 1.4}}$\big)$ is non-decreasing. \end{pro}

{\it Proof}\quad\rm (a)
Firstly, we show that $f_{1}\in L^{1}(\mu)$.
Recall that $\varphi(x)=\nu (0, x)$.
Clearly, for arbitrarily fixed $x_0\in (0, D)$, we have
\footnote{By the integration by parts formula and $\varphi(x)\mu(x,D)\leqslant\delta$, we have $\int_{x_0}^D\!\!\! \sqrt{\varphi}\,\d\mu\leqslant {2\delta}\big/{\sqrt{\varphi(x_{0})}}<\infty.$}
$$\mu\big(\sqrt{\varphi}\,\big)
=\int_{0}^{x_{0}}\!\!\!\sqrt{\varphi}\,\d\mu
+\int_{x_0}^D\!\!\! \sqrt{\varphi}\,\d\mu
\leqslant \int_{0}^{x_{0}}\!\!\!\sqrt{\varphi}\,\d\mu
+\frac{2\delta}{\sqrt{\varphi(x_{0})}}<\infty.$$
Hence $\sqrt{\varphi}\in L^{1}(\mu)$.

(b) Define two sequences $\{h_n\}$ and $\big\{\tilde f_n\big\}$ by
the same recurrence $h_n=h_{n-1}I\!I(h_{n-1})$ but
different initial condition:
\begin{equation*}h_0=1,\qquad
{\tilde f}_{1}=f_{1}=\sqrt{\varphi}. \end{equation*}
We now study $\big\{\tilde f_n\big\}$ first.
From \rf{r5}{Theorem 1.2\,(1)}, we have known that ${\tilde f}_{2}\leqslant 4\delta{\tilde f}_{1}$.
Assume that ${\tilde f}_{n-1}\leqslant (4\delta)^{n-2}{\tilde f}_{1} $ for some $n \geqslant 3$.
Then
$${\tilde f}_{n}\!=\!\int_{0}^{\cdot}\nu(\d y)\int_{y}^{D}{\tilde f}_{n-1}\d\mu
\leqslant\!(4\delta)^{n-2}\int_{0}^{\cdot}\nu(\d y)\int_{y}^{D}{\tilde f}_{1}\d\mu
\!=\!(4\delta)^{n-2}{\tilde f}_{2}
\leqslant (4\delta)^{n-1}{\tilde f}_{1}.$$
By induction, this estimate holds for  $n\geqslant 2$.
Hence ${\tilde f}_{n}\in L^1(\mu)$ for  $n\geqslant 1$
by (a).

Next, we study the sequence $\{h_n\}$. Fix $x_{0}\in(0,D)$. For $x>x_{0}$, we have
\footnote{Since $\sqrt{\varphi}\,\uparrow$, we have $\int_{y}^{D}\sqrt{\varphi}\,\d \mu\geqslant\sqrt{\varphi}\mu(y,D)$. So $$\aligned \int_{x_{0}}^{x}\nu(\d y)\mu(y,D)\leqslant\int_{x_{0}}^{x}\frac{1}{\sqrt{\varphi(y)}}\int_{y}^{D}\sqrt{\varphi}\,\d \mu\,\nu(\d  y)\leqslant\frac{1}{\sqrt{\varphi(x_0)}}\int_{x_0}^{x}\int_y^D\!\!\sqrt{\varphi}\,
\d\mu\,\nu(\d y)\leqslant\frac{{\tilde f}_{2}(x)}{\sqrt{\varphi(x_{0})}}\endaligned$$}
 $$\begin{aligned}
h_{1}(x)&=h_{1}(x_{0})+\int_{x_{0}}^{x}\!\nu(\d  y)\,\mu(y,\,D)\\&
\leqslant h_{1}(x_{0})+\frac{1}{\sqrt{\varphi(x_{0})}}{\tilde f}_{2}(x)
\\&\leqslant h_{1}(x_{0})+\frac{4\delta}{\sqrt{\varphi(x_{0})}}{\tilde f}_{1}(x).
\end{aligned}$$
By induction, it is not difficult to verify that
$$\aligned h_{n}(x)\leqslant \dsum_{k=1}^{n}\frac{{(4\delta)}^{k}}{\sqrt{\varphi(x_{0})}}h_{1}^{n-k}(x_0){\tilde f}_{1}(x)+h_{1}^{n}(x_{0}).\endaligned$$
Hence $h_{n}\in L^{1}(\mu)$ for $n\geqslant1 $.

(c) Now we look for the relationship between $f_{n}$ and ${\tilde f}_{n}$. We begin with
$$\begin{aligned}
f_{1}&={\tilde f}_{1}=\sqrt{\varphi},\\
f_{2}&=\int_{0}^{\cdot}\nu (\d y)\int_{y}^{D}\bar{f}_{1}\d\mu={\tilde f}_{2}-\pi(f_{1})\,h_{1}(x).\end{aligned}$$
By induction, we have in general
$$f_{n}={\tilde f}_{n}-\dsum_{k=1}^{n-1}h_{n-k}\,\pi(f_{k})\qquad n\geqslant 2.$$
Thus $f_{n}\in L^{1}(\mu)$ for every $n\geqslant 1$ by (b).

(d) We now come to the central part of the proof: showing the monotonicity of $\eta_{n}$. By definition of $f_{n}$, we have
\begin{equation}\label{f1.1.18}
\eta_{n}=\sup\limits_{x\in(0,D)}\frac{e^{-C(x)}}{{\bar{f}_{n}}'(x)}\int_{x}^{D}\bar{f}_{n}\,\d\mu
=\sup\limits_{x\in(0,D)}\bigg({\int_{x}^{D}\bar{f}_{n}\,\d \mu}\bigg)\bigg({\int_{x}^{D}\bar{f}_{n-1}\d\mu}\bigg)^{-1}.
\end{equation}
Thus, $\eta_{n}\leqslant\eta_{n-1}$ iff
$$\int_{x}^{D}\big(\bar{f}_{n}-\eta_{n-1}\bar{f}_{n-1}\big)\d\mu
\leqslant 0,\qquad x\in[0, D).$$ That is
\begin{equation*}
\int_{x}^{D}(f_{n}-\eta_{n-1}f_{n-1})\d\mu
\leqslant\big(\pi(f_{n})-\eta_{n-1}\pi(f_{n-1})\big)\,\mu(x,D),\end{equation*}
or equivalently,
\begin{equation}
S(x):=\frac{1}{\mu(x,D)}\int_{x}^{D}(\eta_{n-1}f_{n-1}-f_{n})\d\mu
\geqslant \eta_{n-1}\pi(f_{n-1})-\pi(f_{n})=S(0).
\label{f1.1.10}\end{equation}
This is our key observation and leads to the study on the
monotonicity of $S$.

(e) In view of \eqref{f1.1.10}, we have reduced our proof to showing non-decreasing property of  $S$.
For this, it is enough to show that
\begin{equation*}\mu(y, D)\int_{x}^{D}(\eta_{n-1}f_{n-1}-f_{n})\,\d\mu\leqslant
\mu(x,D)\int_{y}^{D}(\eta_{n-1}f_{n-1}-f_{n})\,\d\mu
\end{equation*}
for any $x, y\in [0, D)$ with $ x<y$. By separating $f_n$ and
$f_{n-1}$, the last inequality is equivalent
to the following one:
\begin{equation}\label{f1.1.11}\eta_{n-1}\!{\int_{y}^{D}}\!\mu(\d u){\int_{x}^{y}}\!\big(f_{n-1}(t)-f_{n-1}(u)\big)\,\mu(\d t)\!\leqslant\!\!
{\int_{y}^{D}}\!\mu(\d u){\int_{x}^{y}}\big(f_{n}(t) -f_{n}(u)\big)\,\mu(\d t).\;\;
\end{equation}
To see this, it suffices to check that
$$f_{n}(u)-f_{n}(t)\leqslant\eta_{n-1}\big(f_{n-1}(u)-f_{n-1}(t)\big),
\qquad u\geqslant t.$$
To check the last inequality, consider $n\geqslant 3$ first. Then
$$\begin{aligned}
f_{n}(u)-f_{n}(t)&=\int_{t}^{u}\nu(\d y)\int_{y}^{D}\bar{f}_{n-1}\d\mu\quad(\text{by definition of $f_n$})\\
&\leqslant\eta_{n-1}\int_{t}^{u}\nu(\d y)\int_{y}^{D}\bar{f}_{n-2}\,\d\mu\quad(\text{by \eqref{f1.1.18}})\\
& =\eta_{n-1}\big(f_{n-1}(u)-f_{n-1}(t)\big)\quad(\text{by definition of $f_{n-1}$}),\qquad u\geqslant t.\end{aligned}$$
It remains to check the required inequality for $n=2$.
 By definition of  $\eta_{1}$, we have
$$\aligned
\frac{e^{-C(y)}}{ \bar{f}_{1}'(y)}\int_{y}^{D}\bar{f}_{1}\,\d\mu=
I\big(\bar{f}_{1}\big)(y)
 \leqslant \eta_{1}.\endaligned$$
It follows that
\begin{equation*}f_{2}(u)-f_{2}(t)=\int_{t}^{u}\nu(\d y)\int_{y}^{D}\bar{f}_{1}\,\d\mu
\leqslant\eta_{1}\int_{t}^{u}f_{1}'(y)\d y\leqslant\eta_{1}\big(f_{1}(u)-f_{1}(t)\big),\quad u\geqslant t.
\end{equation*}
We have thus completed the proof of the monotonicity of
$\{\eta_n\}$ in the continuous context.
$\qquad \Box$

The monotonicity of $\{\eta_n\}$ means we can theoretically improve our lower
estimates of $\lambda_1$ step by step. There is a similar result for the upper estimates but omitted here. It is regretted that the converges of $\{\eta_n^{-1}\}$ to $\lambda_1$ (as $n\to\infty$) remains open. All examples we have ever computed  support the convergence.

\bigskip
{\small
\noindent{\bf Acknowledgements}  The work was supported in part by NNSFC (Grant No. 11131003), SRFDP (Grant No. 20100003110005), the ``985'' project from the Ministry of Education in China and the Fundamental Research Funds for the Central Universities.
}\noindent \\[4mm]

\noindent{\bbb{References}}



\newpage
\subsection*{Appendix A\quad Complement of the proofs in Section 3}
\subsection*{A.1\quad Complementary proof of the two circle arguments: lower estimates}
Let us review the circle arguments for the lower estimates first.
$$\aligned
\lambda_{0}\geqslant\tilde{\lambda}_{0}&
\geqslant\sup\limits_{f\in\mathscr{F}_{I}}\inf\limits_{x\in(0,D)}I(f)(x)^{-1}\\
&=\sup\limits_{f\in\mathscr{F}_{I}}\inf\limits_{x\in(0, D)}I\!I(f)(x)^{-1}=\sup\limits_{f\in\mathscr{F}_{I\!I}}\inf\limits_{x\in(0, D)}I\!I(f)(x)^{-1}\\&\geqslant\sup\limits_{h\in\mathscr{H}}\inf\limits_{x\in(0, D)}R(h)(x)\geqslant\lambda_{0}, \endaligned$$
where $$\tilde{\lambda}_{0}:=\inf\big\{D(f):\mu\big(f^{2}\big)=1,\, f\in{\scr C}^{1}(0, D)\cap{\scr C}[0,D],\, f'(0)=0,\; \text{and}\,f(D)=0\big\}.$$
We prove the circle arguments through the following (a)-(e) steps.

(a) Prove that $\lambda_{0}\geqslant\tilde{\lambda}_{0}\geqslant\sup_{f\in\mathscr{F}_{I\!I}}\inf_{x\in(0,
D)}I\!I(f)(x)^{-1}$.

The first assertion is obvious by definitions of $\lambda_{0}$ and $\widetilde{\lambda}_{0}$. The second one is proved in the main text.

(b)  Prove that $$\sup\limits_{f\in\mathscr{F}_{I}}\inf\limits_{x\in(0, D)}I(f)(x)^{-1}=\sup\limits_{f\in\mathscr{F}_{I}}\inf\limits_{x\in(0,D)}I\!I(f)(x)^{-1}=\sup\limits_{f\in\mathscr{F}_{I\!I}}\inf\limits_{x\in(0, D)}I\!I(f)(x)^{-1}.$$

 For $f\in\mathscr{F}_{1}$,  without loss of generality, assume that $\sup_{x\in(0,D)}I(f)(x)<\infty$. By using Cauchy's mean value theorem, we have
$$\aligned
\sup\limits_{x\in(0, D)}I\!I(f)(x)&=\sup\limits_{x\in(0, D)}\frac{1}{f(x)}\int_{x}^{D}\nu(\d s)\int_{0}^{s}f\d\mu\\
&\leqslant\sup\limits_{x\in(0, D)}-\bigg({\int_{x}^{D}e^{-C(s)}\d s\int_{0}^{s}f\d \mu}\bigg)\bigg({\int_{x}^{D}f'(t)\d t}\bigg)^{-1}\\
&\leqslant \sup\limits_{t\in(0, D)}-\frac{e^{-C(t)}}{f'(t)}\int_{0}^{t}f\d\mu\\
&= \sup\limits_{x\in(0, D)}I(f)(x)<\infty.
\endaligned$$
Making infimum with respect to $f\in\mathscr{F}_{1}$, we have $$\inf\limits_{f\in\mathscr{F}_{I}}\sup\limits_{x\in(0,D)}I\!I(f)(x)\leqslant \inf\limits_{f\in\mathscr{F}_{I}}\sup\limits_{x\in(0, D)}I(f)(x).$$
 Since $\mathscr{F}_{I}\subset\mathscr{F}_{I\!I}$, the left-hand side is bounded below by  $\inf_{f\in\mathscr{F}_{I\!I}}\!\!\sup_{x\in(0, D)}\!\!I\!I(f)\!(x)$. Hence
$$\sup\limits_{f\in\mathscr{F}_{I}}\inf\limits_{x\in(0,D)}I(f)(x)^{-1}\leqslant \sup\limits_{f\in\mathscr{F}_{I}}\inf\limits_{x\in(0, D)}I\!I(f)(x)^{-1}\leqslant
\sup\limits_{f\in\mathscr{F}_{I\!I}}\inf\limits_{x\in(0, D)}I\!I(f)(x)^{-1}.$$

To obtain the equality signs, it suffices to show
$$\sup\limits_{f\in\mathscr{F}_{I\!I}}\inf\limits_{x\in(0, D)}I\!I(f)(x)^{-1}
\leqslant\sup\limits_{f\in\mathscr{F}_{I}}\inf\limits_{x\in(0, D)}I(f)(x)^{-1}.$$
To do so, let $f\in \mathscr{F}_{I\!I}$. Without loss of generality,
assume that $\inf_{x\in(0, D)}I\!I(f)\\(x)^{-1}>0$.
Then $f\in {\scr {C}}[0, D]$ and $f>0$.  Put
$$g(x)=fI\!I(f)(x)=\int_{x}^{D}\nu(\d s)\int_{0}^{s}f\d\mu. $$
Then $g\in\mathscr{F}_{I}$ and $$-g'(s)e^{C(s)}=\int_{0}^{s}f\d\mu\geqslant\int_{0}^{s}g\d\mu\inf\limits_{x\in(0,
D)}\frac{f(x)}{g(x)}\qquad \text{for  }s\in(0,D).$$
That is $I(g)(s)^{-1}\geqslant\inf_{x\in(0, D)}I\!I(f)(x)^{-1}$.
Making infimum with respect to $s\in(0, D)$, we obtain
$$\inf_{s\in (0, D)}I(g)(s)^{-1}\geqslant\inf_{x\in(0, D)}I\!I(f)(x)^{-1}.$$
The assertion now follows by making supremum  with respect to $g\in\mathscr{F}_{I}$  on the both sides of the inequality first and then with respect to $f\in\mathscr{F}_{I\!I}$.

 A different way to prove the equalities here and in (a), without using the continuity of $a$ and $b$, is to show that
 $$\aligned \sup_{f\in\mathscr{F}_{I}}\inf_{x\in (0, D)}I(f)(x)^{-1}\geqslant\lambda_{0}.\endaligned$$

 By the comments below Lemma \ref{L1.1.3} and \eqref{ff1}, we have seen that $\lambda_{0}=I(g)(x)^{-1}$ for $x\in(0, D)$.
In view of Lemmas \ref{L1.1.2} and \ref{L1.1.3}, it follows that $g\in\mathscr{F}_{I}$ and
$$\lambda_{0}=\inf\limits_{x\in(0, D)}I(g)(x)^{-1}\leqslant
\sup\limits_{f\in\mathscr{F}_{I}}\inf\limits_{x\in(0, D)}I(f)(x)^{-1}.$$

{\it{In the following two steps, assume that $a, b\in {\scr C}[0,D]$}.}
\medskip

(c) We show that $\sup_{f\in\mathscr{F}_{I\!I}}\inf_{x\in(0, D)}I\!I(f)(x)^{-1} \geqslant
\sup_{h\in\mathscr{H}}\inf_{x\in(0, D)}R(h)(x)$.

To this end, recall that for each $h\in\mathscr{H}$ with $h={g'}/{g}$ (see Remark \ref{rem1}), we have
$$R(h)=-(ah^{2}+bh+ah')=-\frac{L g}{g}.$$

Before moving further, we prove that if $R(h)>0$ for a positive $g$ with $g'(0)=0$ and $h=g'/g$, then $g$  must be strictly decreasing. In fact, we have $R(h)=-{L g}/{g}>0$. Let $f=gR(h)>0$. Then
$Lg=-f$. Moreover, $g'(x)=-e^{C(x)}\int_{0}^{x}f\d\mu$ since $g'(0)=0$. So $g'<0$ on $(0,D)$.

Now, we return to our main assertion. It suffices to show that $$\aligned\inf_{x\in(0, D)}R(h)(x)\leqslant\sup_{f\in\mathscr{F}_{I\!I}}\inf_{x\in(0,D)}I\!I(f)(x)^{-1}, \qquad h\in\mathscr{H}.\endaligned$$ Without loss of generality, assume that $\inf_{x\in(0, D)}R(h)(x)>0$. Then
$$R(h)(x)>0\qquad\text{ for } x\in(0, D).$$
 Let $f=-(ag''+bg')=gR(h)\;(g\,\text{is the function given above})$. Then $Lg=-f$, $f>0$ and
$f\in{\scr C}[0, D]$  since $a,\,b\in {\scr C}[0,D]$. Since $Lg=-f$ and $g'(0)=0$,  we obtain
$$g(x)-g(D)=\int_{x}^{D}\nu(\d s)\int_{0}^{s}f\d\mu=f(x)I\!I(f)(x)$$
by \eqref{fnew2}. That is $g(x)\geqslant f(x)I\!I(f)(x)$  since $g(D)\geqslant0$. So
$$R(h)(x)^{-1}=\frac{g(x)}{f(x)}\geqslant I\!I(f)(x)\qquad\text{for }x\in(0,D).$$
Furthermore,
$$\inf_{x\in(0, D)}R(h)(x)\leqslant\inf_{x\in(0, D)}I\!I(f)(x)^{-1}
\leqslant\sup_{f\in\mathscr{F}_{I\!I}}\inf_{x\in(0,D)}I\!I(f)(x)^{-1},$$ and
the assertion follows since $h\in\mathscr{H}$ is arbitrary.
\medskip

(d) Prove that $\sup_{h\in\mathscr{H}}\inf_{x\in(0,D)}R(h)(x)\geqslant\lambda_{0}$.
\medskip

Let $f\in L^{1}(\mu),$ $g=fI\!I(f)$, and $\bar{h}={g'}/{g}$ on $[0,D)$. Then $$\bar{h}\in\mathscr{H},\qquad L g=-f,\qquad\text{and}\qquad R(h)=\frac{-L g}{g}=\frac{f}{g}>0.$$ Thus, $\sup_{h\in\mathscr{H}}\inf_{x\in(0,D)}R(h)(x)\geqslant0$. Without loss of generality, assume that $\lambda_{0}>0$. Since $a,\,b\in {\scr C}[0,D]$, by Lemma \ref{L1.1.0} (2), there exists an eigenfunction $g$ such that $Lg=-\lambda_{0}g$. Furthermore,
$$\aligned
g'(0)=0,\quad g|_{(0,D)}>0,\quad g'|_{(0,D)}<0\quad \text{and}\quad g\in{\scr C}^{1}(0,D)\cap{\scr C}[0,D].
\endaligned$$
Let $h={g'}/{g}$. Then $h\in{\scr C}^{1}(0, D)\cap {\scr C}[0,D]$,
$h(0)=0$, $h\in\mathscr{H}$, and $$R(h)(x)=-\frac{Lg(x)}{g(x)}=\lambda_{0}\qquad\text{for } x\in(0, D).$$ So the assertion follows immediately.

(e) We now prove that the supremum in the first circle arguments can be attained. The case that $\lambda={0}$ is easier since
$$\aligned
0=\lambda_{0}\geqslant\inf\limits_{x\in(0,D)}I\!I(f)(x)^{-1}\geqslant0\ \ \text{and}\ \ \  0=\lambda_{0}\geqslant\inf\limits_{x\in(0,D)}I(f)(x)^{-1}\geqslant0
\endaligned$$
for every $f$ in their corresponding domains, as an application of the first circle arguments. Similarly, the conclusion holds for $R$ as seen from proof (d): noting in the degenerated case that $\nu(0,D)=\infty$, we have $\lambda_{0}=0$ (which is a simple consequence of definition \eqref{f1.1}, see also the proof of Theorem \ref{t11.2.1} given below) and then $h=0$ since the eigenfunction is constant in the case.

Next, we consider the case that $\lambda_{0}>0$. Let $g$ be its eigenfunction. For $R$ the supremum  is attained at $h={g'}/{g}$ as seen from the last paragraph of proof (d). For the operator $I$ and $I\!I$, we have already seen that $I(g)\equiv I\!I(g)\equiv\lambda_{0}^{-1}$ according to Lemma \ref{L1.1.4} and the remarks below Lemma \ref{L1.1.3}.\qquad$\Box$\\[4mm]

\subsection*{A.2\quad Complementary proof of the two circle arguments: upper estimates}
For the upper estimation of $\lambda_{0},$  we review and show the circle arguments in the following.
$$\aligned
\lambda_{0}&\leqslant\inf\limits_{f\in\mathscr{\widetilde{F}}_{I\!I}\bigcup\mathscr{\widetilde{F}}_{I\!I}'}
\sup\limits_{x\in \supp(f)}I\!I(f)(x)^{-1}
=\inf\limits_{f\in\mathscr{\widetilde{F}}_{I\!I}}\sup\limits_{x\in
\supp(f)}I\!I(f)(x)^{-1}
\\&= \inf\limits_{f\in\mathscr{\widetilde{F}}_{I}}\sup\limits_{x\in
\supp(f)}I\!I(f)(x)^{-1}=
\inf\limits_{f\in\mathscr{\widetilde{F}}_{I}}\sup\limits_{x\in(0, D)}I(f)(x)^{-1}
\\&=\inf\limits_{f\in\mathscr{\widetilde{F}'}_{I}}\sup\limits_{x\in(0,
D)}I(f)(x)^{-1}\\&\leqslant \inf\limits_{h\in\mathscr{\widetilde{H}}}\sup\limits_{x\in(0,
D)}R(h)(x)\leqslant\lambda_{0}.
\endaligned$$

(f) The assertion that
$\lambda_{0}\leqslant\inf_{f\in\mathscr{\widetilde{F}}_{I\!I}\bigcup\mathscr{\widetilde{F}}_{I\!I}'}
\sup_{x\in \rm supp(f)}I\!I(f)(x)^{-1}$ is proved in our main text.
\medskip

(g) Prove that $$\inf\limits_{f\in\mathscr{\widetilde{F}}_{I\!I}}\sup\limits_{x\in
\rm supp(f)}I\!I(f)(x)^{-1}=
\inf\limits_{f\in\mathscr{\widetilde{F}}_{I}}\sup\limits_{x\in
\rm supp(f)}I\!I(f)(x)^{-1}=
\inf\limits_{f\in\mathscr{\widetilde{F}}_{I}}\sup\limits_{x\in(0,
D)}I(f)(x)^{-1}. $$\rm

Let $f\in\mathscr{\widetilde{F}}_{I}$, there exists $x_{0},\ x_{1}\in [0,D)$ such that $f=f(\cdot\vee x_{0})\mathbbold{1}_{[0,x_{1})}\in {\scr C}^{1}(x_{0},x_{1})\cap {\scr C}[x_{0},x_{1}]$.  By Cauchy's mean value theorem, we have
$$\aligned
\inf\limits_{x\in \rm supp(f)}I\!I(f)(x)&=\inf\limits_{x\in [x_{0},x_{1})}\frac{1}{f(x)}\int_{x}^{x_{1}}e^{-C(t)}\int_{0}^{t}f\d\mu\d t\\
&\geqslant\inf\limits_{x\in [x_{0},x_{1})}\frac{1}{-f'(x)}e^{-C(x)}\int_{0}^{x}f\d\mu\\
&=\inf\limits_{x\in (0,D)}I(f)(x).
\endaligned$$
So the assertion that $$
\inf\limits_{f\in\mathscr{\widetilde{F}}_{I\!I}}\sup\limits_{x\in
\rm supp(f)}I\!I(f)(x)^{-1}\leqslant
\inf\limits_{f\in\mathscr{\widetilde{F}}_{I}}\sup\limits_{x\in
\rm supp(f)}I\!I(f)(x)^{-1}\leqslant
\inf\limits_{f\in\mathscr{\widetilde{F}}_{I}}\sup\limits_{x\in (0,
D)}I(f)(x)^{-1}$$ follows by $\widetilde{\mathscr{F}}_{I}\subset\widetilde{\mathscr{F}}_{I\!I}$.

There are two choices to prove the equalities.
The first choice is  proving the assertion that $$\inf\limits_{f\in\mathscr{\widetilde{F}}_{I}}\sup\limits_{x\in
(0,D)}I(f)(x)^{-1}\leqslant\inf\limits_{f\in\mathscr{\widetilde{F}}_{I\!I}}\sup\limits_{x\in
\rm supp(f)}I\!I(f)(x)^{-1}.$$
For  $f\in\widetilde{\mathscr{F}}_{I\!I},\ \exists
x_{0}\in(0,D)$ such that $f=f \mathbbold{1}_{[0, x_{0})}$ and $f\in{\scr C}[0,x_{0}]$.  Let $g=fI\!I(f) \mathbbold{1}_{\rm supp(f)}$.
Then
\begin{equation*}g(x)=\int_{x}^{x_{0}}\nu(\d s)\int_{0}^{s}f\d\mu\mathbbold{1}_{[0, x_{0})}(x)
,\qquad g\in\widetilde{\mathscr{F}}_{1}'\subseteq\widetilde{\mathscr{F}}_{I},\
\end{equation*} and $Lg=-f$ on $[0, x_{0})$ by simple
calculation. Since $g'(0)=0$, replacing $[0,D]$ with $[0, x_{0})$ in \eqref{fnew1},
we have
\begin{equation*}
-e^{C(x)}g'(x)=\int_{0}^{x}f\d\mu\leqslant\int_{0}^{x}g\d\mu\sup\limits_{t\in(0,
x_{0})}\frac{f(t)}{g(t)},\qquad x<x_0.\ \end{equation*} Hence,
\begin{equation*}-e^{C(x)}g'(x)\bigg({\int_{0}^{x}g\d\mu}\bigg)^{-1}\leqslant\sup\limits_{x\in
\rm supp(f)}I\!I(f)(x)^{-1},\qquad x<x_0. \end{equation*} Making supremum with respect
to $x\in(0, x_{0})$,\, we have
\begin{equation*}\sup\limits_{x\in (0, D)}I(g)(x)^{-1}=\sup\limits_{x\in (0,
x_{0})}I(g)(x)^{-1}\leqslant\sup\limits_{x\in \rm supp(f)}I\!I(f)(x)^{-1}.
\end{equation*}The assertion now follows by making infimum with respect to
$g\in\widetilde{\mathscr{F}}_{I}$ first,\ then with respect to
$f\in\widetilde{\mathscr{F}}_{I\!I}.$
\medskip

The second method for the identity is  making a small circle below.

  Since
$$\aligned
\lambda_{0}&\leqslant\inf\limits_{f\in\mathscr{\widetilde{F}}_{I\!I}}\sup\limits_{x\in
\rm supp(f)}I\!I(f)(x)^{-1}\leqslant
\inf\limits_{f\in\mathscr{\widetilde{F}}_{I}}\sup\limits_{x\in
\rm supp(f)}I\!I(f)(x)^{-1}\\&\leqslant
\inf\limits_{f\in\mathscr{\widetilde{F}}_{I}}\sup\limits_{x\in (0,
D)}I(f)(x)^{-1}\leqslant\inf\limits_{f\in\mathscr{\widetilde{F}}'_{I}}\sup\limits_{x\in
(0, D)}I(f)(x)^{-1},
\endaligned$$ it suffices to show
$\inf\limits_{f\in\mathscr{\widetilde{F}}'_{I}}\sup\limits_{x\in (0,
D)}I(f)(x)^{-1}\leqslant\lambda_{0}. $

To see this, we  introduce an approximating procedure.
Recall that
$$\aligned&\lambda_{0}^{(0,p)}=\inf
\big\{D(f):\mu\big(f^{2}\big)=1,\ f\in\mathscr{C}^{1}(0,
p)\cap\mathscr{C}[0,p],\,f|_{[p,D)}=0\big\}.\endaligned$$
 Let $ p_{n}\in(0, D),\ p_{n}\uparrow D$. Then
$\lambda_{0}^{(0,p_{n})}\downarrow\lambda_{0}$ by  Lemma \ref{L1.1.1}\,(1), where
$\lambda_{0}^{(0,p_{n})}$ is the first eigenvalue of the Dirichlet form $(D,\mathscr{D}(D))$ restricted to
$(0, p_{n})$ with ND-boundaries. Now, let $g$ be the eigenfunction of $\lambda_{0}^{(0,p_{n})}>0.$ Extend $g$ to the whole space by setting $g=g\mathbbold{1}_{[0,p_{n})}.$ By using Lemmas \ref{L1.1.2} and \ref{L1.1.3}, it follows that $g\in\widetilde{\mathscr{F}}_{I}'$. Furthermore,
\begin{equation*}\lambda_{0}^{(0,p_{n})}=\sup\limits_{x\in(0, p_{n})}I(g)(x)^{-1}=
\sup\limits_{x\in(0, D)}I(g)(x)^{-1}\geqslant
\inf\limits_{g\in\widetilde{\mathscr{F}}_{I}'}\sup\limits_{x\in(0,
D)}I(g)(x)^{-1}. \end{equation*} The assertion now follows by letting $n\rightarrow\infty$ because of $\tilde{\lambda}_{0}=\lambda_{0}$.
\medskip

{\it In the following two steps, we assume that $a,b\in\mathscr{C}[0,D]$.}
\medskip

(h) Prove that $\inf_{f\in\mathscr{\widetilde{F}}_{I}}\sup_{x\in(0,
D)}I\!I(f)(x)^{-1}\leqslant
\inf_{h\in\mathscr{\widetilde{H}}}\sup_{x\in(0, D)}R(h)(x)$.
\medskip \rm

Firstly, for
$h\in\widetilde{\mathscr{H}},$
$\exists x_{0}\in(0, D)$ such that  $R(h)>0$ on $(0, x_{0})$. We use $R(g)|_{[0,x_{0})}$ instead  of $R(h)|_{[0,x_{0})}$ as in (c).  Hence,
$$R(g)(x)=\begin{cases}
 -({L g}/{g})(x),& x<x_{0};\\0,& x\geqslant
x_{0}.
\end{cases}$$
Secondly, we turn to the main assertion.

 Let $f=gR(g)$. Then $$f=f \mathbbold{1}_{[0, x_{0})},\qquad\qquad
 L g=-f\quad\text{on}\; (0, x_{0}),$$ and $f\in\widetilde{\mathscr{F}}_{I\!I}$ since $a,\, b\in{\scr C}[0,D]$.
Noting that $g'(0)=0$ and $g(x_{0})=0$, by \eqref{fnew2}, we have\begin{equation*}
g(y)=\int_{y}^{x_{0}}\nu(\d x)\int_{0}^{x}f\d\mu=f(y)I\!I(f)(y)\qquad\text{for }\quad y<x_{0}.
\end{equation*}
So $$R(h)(y)=\frac{f(y)}{g(y)}=I\!I(f)(y)^{-1}\qquad\text{for }\quad
y< x_{0}.$$ Making supremum with respect to $y\in(0, x_{0})$,\ we
have\begin{equation*}\sup\limits_{y\in(0, x_{0})}R(h)(y)=\sup\limits_{y\in(0,
x_{0})}I\!I(f)(y)^{-1},
\end{equation*} and  the assertion follows immediately by making infimum with respect to $h\in\widetilde{\mathscr{H}}$ first and then making infimum with respect to $f\in\widetilde{\mathscr{F}}_{I\!I}$.
\medskip

(i) Prove that $\inf_{h\in\mathscr{\widetilde{H}}}\sup_{x\in(0,D)}R(h)(x)\leqslant\lambda_{0}.$
\medskip

When $D<\infty$, since  $a,b\in{\scr C}[0,D]$, there is an eigenfunction $g$ satisfying
$$h:=\frac{g'}{g}\in\widetilde{\mathscr{H}},\qquad\quad R(h)=-\frac{L g}{g}=\lambda_{0}.$$
Indeed, since $a,b\in{\scr C}[0,D],$ we have $g\in{\scr C}^{2}[0, D]$, $g'(0)=0$, $g(D)=0$, and $g'<0$ on $(0, D)$. Hence,\, ${h}(0)=0,\;{h}(D)=0,\;{h}<0\;\text{on}\;(0, D),\;\text{and}\;
{h}\in\mathscr{C}^{1}(0, D)\cap{\scr C}[0,D]$. Moreover, $R({h})(x)=-\big({L g}/{g}\big)(x)=\lambda_{0}>0$. So the assertion holds for $D<\infty$.

When $D=\infty$, let $p_{n}\uparrow \infty$. For fixed $p_{n}$, as the last part of \ (g), denote by $g$ the eigenfunction of $\lambda_{0}^{(0,p_{n})}>0$, i.e. $$\aligned L g(x)=-\lambda_{0}^{(0,p_{n})}g(x), \qquad x\in(0,p_{n}).\endaligned$$ Since $a,b\in {\scr C}[0,D]$, we have
 $$\aligned g\downarrow\downarrow\;\text{on}\;(0,p_{n}),\quad g'(0)=0,\quad g(p_{n})=0,\quad g\in\mathscr{C}^{2}[0,p_{n}].\endaligned$$
 by Lemmas \ref{L1.1.2}, \ref{L1.1.3} and \ref{L1.1.4}.

 Let $\bar{h}_{n}(x)=-{g'(x)}\mathbbold{1}_{[0,p_{n})}(x)/{g(x)}$.
Then $\bar{h}_{n}(x)\in\mathscr{\widetilde{H}}$ and
$$\aligned
\lambda_{0}^{(0,p_{n})}&=\sup\limits_{x\in(0, p_{n})}R(g)(x)\geqslant
 \inf\limits_{h\in\mathscr{\widetilde{H}},\rm supp(h)=(0, p_{n})}\sup\limits_{x\in(0, p_{n})}R(h)(x)\\
&\geqslant
 \inf\limits_{h\in\mathscr{\widetilde{H}},\rm supp(h)=(0, p_{n})}\sup\limits_{x\in(0, D)}R(h)(x)\\&\geqslant
 \inf\limits_{h\in\mathscr{\widetilde{H}}}\sup\limits_{x\in(0, D)}R(h)(x). \endaligned$$
The assertion now follows by letting $n\rightarrow
 \infty$.$\qquad \Box$

\subsection*{A.3\quad Proof of Proposition $\ref{t11.1.3}$}

The proof consists of the following four parts.

(a) The assertion that $$\sup_{h\in\overline{\mathscr{H}}}\inf_{x\in(0, D)}\overline R(h)(x)\geqslant\sup_{h\in\mathscr{H}_{\ast}}\inf_{x\in(0,D)}\overline R(h)(x)\geqslant\lambda_{0}\geqslant0$$ is proved in our main text.

(b) Prove that $\lambda_{0}=\sup_{h\in\overline{\mathscr{H}}}\inf_{x\in(0, D)}\overline R(h)(x)$ whenever $\mu(0, D)=\infty$.\\
From (a), it suffices to show that
\begin{equation*}\lambda_{0}\geqslant\sup\limits_{h\in\overline{\mathscr{H}}}\inf\limits_{x\in(0, D)}\overline R(h)(x),\end{equation*}
or equivalently $$\aligned\lambda_{0}\geqslant\inf_{x\in(0, D)}\overline R(h)(x)\qquad\text{for every} \;h\in\overline{\mathscr{H}}.\endaligned$$
In view of (a), without loss of generality, assume that $\inf_{x\in(0, D)}\overline R(h)(x)>0$ for a given $h\in\overline{\mathscr{H}}$. Let $f=-(ah'+b h)$. Since $h<0$ on $(0, D)$ and $h(0)=0$,  we have
$$\aligned f'=h \overline R(h)<0\qquad\text{and }\qquad f'(0)=0.\endaligned$$
Thus, $ f\in{\scr C}^{1}(0, D)\cap {\scr C}[0,D]$. It follows that $f\in\mathscr{F}_{I}$ once we show that $f>0$ on $(0,D)$.

For this, fix $ x\in (0, D)$. By integration formula by parts and $h(0)=0$, we obtain
\begin{equation}\label{f1.1.5}\int_{0}^{x}f\d\mu=-e^{C(x)}h(x)>0. \end{equation}
Since $f\;\downarrow\downarrow$, if $f(x_{0})\leqslant0$ for some $x_{0}\in (0,D)$, then
\begin{equation*}\int_{x_{0}}^{x}f\d\mu\leqslant f(x_{0})\mu(x_{0},x),\qquad x>x_0,\end{equation*}
and the right-hand side of the above inequality converges to $-\infty$ as $x\rightarrow D$.

By \eqref{f1.1.5}, we obtain
\begin{equation*}0<\int_{0}^{x}f\d\mu=\int_{0}^{x_{0}}f\d\mu+\int_{x_{0}}^{x}f\d\mu
\rightarrow-\infty\qquad \text{as}\ x\rightarrow D, \end{equation*}
which is a contradiction. So $f>0$ on $(0,D)$.

Because of \eqref{f1.1.5}, we have
\begin{equation*}\overline R(h)(x)=\frac{f'(x)}{h(x)}=\bigg(-\frac{e^{-C(x)}}{f'(x)}\int_{0}^{x}f\d\mu\bigg)^{-1}=I(f)(x)^{-1}. \end{equation*}
Making infimum with respect to $x\in(0, D)$ first and then making  supremum with respect to $f\in \mathscr{F}_{I}$, we obtain the assertion by the variational formulas for lower bounds in Theorem \ref{th1}\,(1).
\medskip

(c) Prove  that $\lambda_{0}=\sup_{h\in\mathscr{H}_{\ast}}\inf_{x\in(0, D)}\overline R(h)(x)$.

It suffices to prove that $$\lambda_{0}\geqslant\sup_{h\in\mathscr{H}_{\ast}}\inf_{x\in(0, D)}\overline R(h)(x).$$ The main body in proof (b) is to prove that the function $f$ defined there is positive, this is automatic due to the definition of  $\mathscr{H}_{\ast}$. So the proof (b) can be applied to $h\in\mathscr{H}_{\ast}$ directly.
Hence, $\lambda_{0}\geqslant\sup_{h\in\mathscr{H}_{\ast}}\inf_{x\in(0, D)}\overline R(h)(x)$ and then the equality $\lambda_{0}=\sup_{h\in\mathscr{H}_{\ast}}\inf_{x\in(0, D)}\overline R(h)(x)$ follows.$\qquad \Box$

\subsection*{ A.4\quad Proof of Lemma $\ref{l00}$}
Define $M(x)=\int_{0}^{x}m(y)\d y$. Using integration by parts formula, we have
$$\aligned
&\int_{0}^{x}m(y)\psi^r(y)\d y=\int_{0}^{x}\psi^r(y)\d M(y)\\
&\qquad=\psi^r(x)M(x)-r\int_{0}^{x}\psi^{r-1}(y)\psi'(y)M(y)\d y\quad (\text{by}\ M(0)=0,\, \psi(0)<\infty)\\
&\qquad\leqslant c\psi^{r-1}(x)-cr\int_{0}^{x}\psi^{r-2}(y)\psi'(y)\d y\quad (\text{since}\, M\psi\leqslant c, \psi'=-n\leqslant 0)\\
&\qquad=c\psi^{r-1}(x)-\frac {cr}{r-1}(\psi^{r-1}(x)-\psi^{r-1}(0))\\
&\qquad=\frac{c}{1-r}\psi^{r-1}(x)-\frac {cr}{1-r}\psi^{r-1}(0)\leqslant\frac{c}{1-r}\psi^{r-1}(x). \qquad\Box
\endaligned$$

\subsection*{A.5\quad Proof of Theorem $\ref{t11.2.1}$}

\rm Firstly, the assertion $\lambda_{0}\geqslant (4{\delta})^{-1}$ is proved in the main text.

Now, we show that $\lambda_{0}\leqslant\delta^{-1}$.
Let $ x_{0}$, $x_{1}\in[0, D)$ with $x_{0}<x_{1}$.
Set $f=\nu(x_{0}\vee\cdot,x_{1})\mathbbold{1}_{[0, x_{1})}$.
Then $f\in\mathscr{\widetilde{F}}_{I}$, $f'=-e^{-C}$ on $(x_{0}, x_{1})$, and
$$I(f)(x)=\begin{cases}
 \int_{x_{0}}^{x}f\d\mu+\nu(x_{0},x_{1})\mu{(0, x_{0})}, &x\in(x_{0}, x_{1});\\ \infty\ (\text{ by convention, $1/0=\infty$}) ,&\text{otherwise}.
\end{cases}$$
Thus, $I(f)(x)$ achieves its minimum at $\mbox{$x=x_0+$}$, and
 $$\inf\limits_{x\in(0, D)}I(f)(x)=\inf\limits_{x\in(x_{0}, x_{1})}I(f)(x)=\nu(x_{0},x_{1})\mu{(0, x_{0})}\to\mu(0, x_{0})\nu(x_0,D)$$
$\text{as}\ x_1\to D$. Hence,
\begin{equation*}\lambda_{0}^{-1}\geqslant\sup\limits_{f\in\mathscr{\widetilde{F}}_{I}}\inf\limits_{x\in(0, D)}I(f)(x) \geqslant\sup\limits_{x_{0}\in(0, D)}\mu{(0, x_{0})}\nu(x_0, D) =\delta. \end{equation*}

Another method to show that $\lambda_{0}\leqslant\delta^{-1}$ is  using the classical variational formula:
$$\lambda_{\ast}^{-1}\geqslant\sup_{x_{0},x_{1}:x_{0}<x_{1}}
\frac{\big\|f^{x_{0},x_{1}}\big\|}{D\big(f^{x_{0},x_{1}}\big)}.$$
By simple calculation, we have
$$\aligned
&\big\|f^{x_{0},x_{1}}\big\|=\int_{0}^{D}{f^{x_{0},x_{1}}}^{2}\d\mu
=\nu^{2}(x_{0},x_{1})\mu(0,x_{0})+\int_{x_{0}}^{x_{1}}\nu^{2}(t,x_{1})\d\mu,\\
&D\big(f^{x_{0},x_{1}}\big)=\nu(x_{0},x_{1}).
\endaligned$$
Thus,
$$\aligned\lambda_{0}^{-1}=\tilde{\lambda}_0^{-1}=\lambda_{\ast}^{-1}\geqslant\!\!\sup_{x_{0},x_{1}:x_{0}<x_{1}}
\frac{\big\|f^{x_{0},x_{1}}\big\|}{D\big(f^{x_{0},x_{1}}\big)}\geqslant\sup_{x_{0},x_{1}:x_{0}<x_{1}}\nu(x_{0},x_{1})\mu(0,x_{0})=\delta\endaligned$$
by Lemma \ref{new} and Proposition \ref{p11.1.1}\,(1).
$\qquad\Box$
\subsection*{A.6\quad Proof of Theorem $\ref{t11.2.2}$ and Corollary $\ref{cor1}$}
 We prove the  assertions through the following six steps.

(a) By Cauchy's  mean value theorem and \eqref{f1.1.16}, we have
\begin{equation*}\delta_{1}=\sup\limits_{x\in(0, D)}I\!I(f_{1})(x)
\leqslant\sup\limits_{x\in(0, D)}I(f_{1})(x)\leqslant4\delta,\end{equation*}
and
\begin{equation*}\delta_{n+1}=\sup\limits_{x\in(0, D)}I\!I(f_{n+1})(x)=\sup\limits_{x\in{(0, D)}}\frac{f_{n+2}(x)}{f_{n+1}(x)}
\leqslant\sup\limits_{x\in(0, D)}\frac{f_{n+1}(x)}{f_{n}(x)}=\delta_{n},
\end{equation*}
which means the monotonicity of  $\delta_{n}$ with respect to $n$.

Notice that
\begin{equation*}f_{1}(x)=\sqrt{\varphi(x)}>0,\qquad f_{1}'(x)=-\frac{e^{-C(x)}}{2\sqrt{\varphi(x)}}. \end{equation*}
We have $f_1\in \scr{C}[0,D]$ and  $f_{1}\in\mathscr{F}_{I\!I}$. Moreover, by  induction, we have
$f_{n}\in\mathscr{F}_{I\!I}$ for $n\geqslant 1$. Therefore,
$$\aligned\lambda_{0}\geqslant\sup\limits_{f\in \mathscr{F}_{I\!I}}\inf\limits_{x\in(0,D)}I\!I(f)(x)^{-1}\geqslant{\delta_{n}^{-1}}
\qquad\text{for}\quad n\geqslant 1.\endaligned$$.

(b) By Cauchy's mean value theorem, we have
$$\begin{aligned}
\inf\limits_{x< x_{1}}I\!I\big(f_{1}^{x_{0},x_{1}}\big)(x)&=\inf\limits_{x\in[x_{0}, x_{1})}\frac{1}{f_{1}^{x_{0},x_{1}}(x)}\int_{x}^{x_{1}}e^{-C(s)}\d s
\int_{0}^{s}f_{1}^{x_{0}, x_{1}}\d\mu\\
&\geqslant\inf\limits_{x\in[x_{0}, x_{1})}\int_{0}^{x}f_{1}^{x_{0},x_{1}}\d\mu =\inf\limits_{x\in[x_{0}, x_{1})}\int_{0}^{x}\int_{t\vee x_{0}}^{x_{1}}\d\nu\mu(\d t)\\
&=\int_{0}^{x_{0}}d\mu\int_{x_{0}}^{x_{1}}\d \nu=\mu(0,x_{0})\nu(x_{0}, x_{1}).\end{aligned}$$
Thus,
$$\aligned
\sup_{x_{0}, x_{1}:x_{0}<x_{1}}\!\!\inf_{x< x_{1}}I\!I\big(f_{1}^{x_{0},x_{1}}\big)(x)\!\!\geqslant\!\!\sup_{x_{0}, x_{1}:x_{0}<x_{1}}\!\!\mu(0,x_{0})\nu(x_{0}, x_{1})
\!=\!\!\sup\limits_{x\in(0, D)}\mu(0,x)\nu(x,D),\endaligned$$
which is just $\delta_{1}'\geqslant\delta$.

Meanwhile, for the same reason, we obtain
\begin{equation*}
\inf\limits_{x<x_{1}}I\!I\big(f_{n+1}^{x_{0}, x_{1}}\big)(x)=\inf\limits_{x< x_{1}}\frac{f_{n+2}^{x_{0}, x_{1}}(x)}{f_{n+1}^{x_{0}, x_{1}}(x)}
\geqslant\inf\limits_{x<{x_{1}}}\frac{f_{n+1}^{x_{0}, x_{1}}(x)}{f_{n}^{x_{0}, x_{1}}(x)}= \inf\limits_{x< x_{1}}I\!I\big(f_{n}^{x_{0}, x_{1}}\big)(x).
\end{equation*}
which implies that
$$\delta_{n+1}'=\sup\limits_{x_{0}, x_{1}:x_{0}< x_{1}}\inf\limits_{x< x_{1}}I\!I\big(f_{n+1}^{x_{0}, x_{1}}\big)(x)\geqslant
\sup\limits_{x_{0}, x_{1}:x_{0}< x_{1}}\inf\limits_{x< x_{1}}I\!I\big(f_{n}^{x_{0},x_{1}}\big)(x)=\delta_{n}',$$
i.e.  $\delta_{n}'$ is increasing in $n$.

(c) Noticing that $$\aligned f_{1}^{x_{0}, x_{1}}=\nu(x_{0}\vee\cdot,x_{1})\mathbbold{1}_{[0,x_{1})}\qquad\text{ and}\qquad \big(f_{1}^{x_{0}, x_{1}}\big)'(x)=-e^{-C(x)}\quad\text{on }(x_0, x_{1}),\endaligned$$ we have
$f_{1}^{x_{0}, x_{1}}\in{\scr C}^{1}(x_0, x_1)\cap {\scr C}[x_0,x_1]$ and further $f_{1}^{x_{0}, x_{1}}\in\mathscr{\widetilde{F}}_{I}\subset\mathscr{\widetilde{F}}_{I\!I}$.
 It is easy to verify that $f_{n}^{x_{0}, x_{1}}\in\mathscr{\widetilde{F}}_{I}\subset\mathscr{\widetilde{F}}_{I\!I}$ for $n\geqslant 1$ by induction. From the variational formula for upper bounds, we obtain the following inequalities below immediately.
$$\lambda_{0}\leqslant \inf\limits_{f\in\mathscr{\widetilde{F}}_{I}}\sup\limits_{x\in \text{supp} (f)}I\!I(f)(x)^{-1}\leqslant
\inf\limits_{x_{0}, x_{1}:x_{0}<x_{1}}\sup\limits_{x< x_{1}}I\!I\big({f_{n}}^{x_{0}, x_{1}}\big)(x)^{-1}\!\!={\delta'_{n}}^{-1},\;\; n\geqslant 1.$$

\medskip
(d)
Since $f_{n}^{x_{0}, x_{1}}$ By  definition of  $\lambda_{\ast}$,
 it is obvious that
 $$\aligned(\bar{\delta}_{n})^{-1}\geqslant \lambda_{\ast}=\tilde{\lambda}_{0}=\lambda_{0}\endaligned$$
 by Lemma \ref{new}  and Proposition \ref{p11.1.1}\,(1).  Next, let $f=f_{n}^{x_{0},x_{1}}$.
  Replacing $fI\!I(f)(\cdot\wedge x_{0})\mathbbold{1}_{\supp(f)}=f_{n+1}^{x_{0},x_{1}}$ with $g$ in \eqref{f1.2}, by definition of $\lambda_{\ast}$, we obtain that $\bar{\delta}_{n+1}\geqslant\delta_{n}'.$
\medskip

(e) The computation of $\delta_{1}$ is simple, Here we only compute $\bar{\delta}_{1} $ and $ \delta_{1}'$ in details.  Firstly, for any $x\in[x_0, x_1)$, we have
$$\begin{aligned}
f_{2}^{x_{0}, x_{1}}(x)&=\int_{x}^{x_{1}}\nu(\d s)\int_{0}^{s}f_{1}^{x_{0}, x_{1}}\d\mu= \int_{0}^{x_{1}}f_{1}^{x_{0}, x_{1}}(t)\mu(\d t)\int_{t\vee x}^{x_{1}}\d \nu\\
&=\int_{0}^{x_{1}}\mu(\d t)\int_{t\vee x_{0}}^{x_{1}}\d\nu\int_{t\vee x}^{x_{1}}\d\nu\\
&=\int_{0}^{x_{0}}\mu(\d t)\int_{x_{0}}^{x_{1}}\d\nu\int_{x}^{x_{1}}\d\nu
        +\int_{x_{0}}^{x}\mu(\d t)\int_{t}^{x_{1}}\d\nu\int_{x}^{x_{1}}\d\nu  \\
&\qquad +\int_{x}^{x_{1}}\mu(\d t)\int_{t}^{x_{1}}\d\nu\int_{t}^{x_{1}}\d\nu\\
&=\nu(x_{0}, x_{1})\nu(x, x_{1})\mu(0, x_{0})+\nu(x,x_{1})\!\int_{x_{0}}^{x}\nu(t, x_{1})\mu(\d t)\\&\qquad+\int_{x}^{x_{1}}\nu(t, x_{1})^2\mu(\d t)\\
&=:\bigg(\nu(x_{0}, x_{1})\mu(0, x_{0})+H_{1}(x)\bigg)\nu(x, x_{1}),\end{aligned}$$
where
\begin{equation*}H_{1}(x)=\int_{x_{0}}^{x}\nu(t, x_{1})\mu(\d t)+\frac{1}{\nu(x, x_{1})}\int_{x}^{x_{1}}\nu^{2}(t, x_{1})\mu(\d t),\qquad x\in [x_0, x_1). \end{equation*}
Noticing that $f_{1}^{x_{0}, x_{1}}(x)=\nu(x, x_{1})$ for every $x\in [x_{0},x_{1})$,
we have
$$\aligned\inf_{x_{0}\!\leqslant\!x< x_{1}}\!\frac{f_{2}^{x_{0},
x_{1}}(x)}{f_{1}^{x_{0}, x_{1}}(x)}\!=\!\nu(x_{0}, x_{1})\mu(0, x_{0})\!+\!\inf_{x_{0}\leqslant x< x_{1}}H_{1}(x)\!=\!\nu(x_{0}, x_{1})\mu(0, x_{0})\!+\!H_{1}(x_0).
\endaligned$$
In the last equality, we have used the fact that $H_{1}$ is non-decreasing on $[x_{0},x_{1})$. Indeed, fix $x, y\in[x_{0}, x_{1})$ with $x< y$. Since $\nu(t, x_{1})$ is decreasing in $t\in(x_0, x_1)$,  we have
\begin{equation*}\frac{1}{\nu(x, x_{1})}\int_{x}^{y}\nu^{2}(t, x_{1})\mu(\d t)\leqslant\int_{x}^{y}\nu(t, x_{1})\mu(\d t)\quad \text{and}\quad \frac{1}{\nu(y, x_{1})}-\frac{1}{\nu(x, x_{1})}>0. \end{equation*}
Moreover,
\begin{eqnarray}\frac{1}{\nu(x, x_{1})}\int_{x}^{y}{\nu^{2}(t, x_{1})\mu(\d t)}\hskip-0.4cm&&\!\!\leqslant\int_{x}^{y}\nu(t, x)\mu(\d t)\nonumber\\&&\quad+
\bigg(\frac{1}{\nu(y, x_{1})}-\frac{1}{\nu(x, x_{1})}\bigg)\!\!\int_{y}^{x_{1}}\!\!{\nu^{2}}(t, x_{1})\mu(\d t),
\label{f1.1.6}\end{eqnarray}
which implies that $H_{1}(x)\leqslant H_{1}(y)$.
 So $H_{1}$ is non-decreasing on $[x_0, x_1)$.
Therefore,
\begin{eqnarray}
\delta_{1}'\hskip-0.6cm&&=\sup\limits_{x_{0},x_{1}:x_{0}<x_{1}}\inf\limits_{x<x_{1}}I\!I(f_{1}^{x_{0}, x_{1}})(x)=
 \sup\limits_{x_{0}, x_{1}:x_{0}<x_{1}} \inf\limits_{x_{0}\leqslant x< x_{1}}\frac{f_{2}^{x_{0},x_{1}}(x)}{f_{1}^{x_{0},
x_{1}}(x)}\nonumber\\
&&= \sup\limits_{x_{0}, x_{1}:x_{0}<x_{1}}\bigg(\nu(x_{0}, x_{1})\mu(0, x_{0})+\inf\limits_{x_{0}\leqslant x<x_{1}}H_{1}(x)\bigg)
\nonumber\\
&&=\sup\limits_{x_{0}, x_{1}:x_{0}<x_{1}}\bigg(\nu(x_{0}, x_{1})\mu(0, x_{0})+\frac{1}{\nu(x_{0},x_{1})}\int_{x_{0}}^{x_{1}}\nu^{2}(t, x_{1})\mu(\d t)\bigg)\label{f1.1.7}\\
&&=\sup\limits_{x_0\in(0, D)}\bigg(\nu(x_{0},D)\mu(0, x_{0})+\frac{1}{\nu(x_{0},
D)}\int_{x_{0}}^{D}\nu^{2}(t,D)\mu(\d t)\bigg). \label{f1.1.17}
\end{eqnarray}
In \eqref{f1.1.17}, we have used the fact that
$$H_{2}(x):=\nu(x_{0}, x)\mu(0, x_{0})+\frac{1}{\nu(x_{0},x)}\int_{x_{0}}^{x}\nu^{2}(t, x)\mu(\d t),\qquad x>x_0$$ is non-decreasing in $x$. In fact,
for  $x_{0}<x< y$,  $H_{2}(x)\leqslant H_{2}(y)$ if and only if
\begin{equation}\nu(x,y)\mu(0, x_{0})+\frac{1}{\nu(x_0,y)}\int_{x}^{y}\nu^{2}(t,y)\mu(\d t)+
\int_{x_{0}}^{x}\bigg(\frac{\nu^{2}(t,y)}{\nu(x_{0},y)}-\frac{\nu^{2}(t, x)}{\nu(x_{0},x)}\bigg)\,\mu(\d t)\geqslant0.
\label{f1.1.8}\end{equation}
For $t\in[x_{0},x]$, we have
$$\frac{\nu^{2}(t, y)}{\nu^{2}(t, x)}=\bigg(1+\frac{\nu(x, y)}{\nu(t,x)}\bigg)^{2}\geqslant1+\frac{\nu(x, y)}{\nu(t, x)}\geqslant1+\frac{\nu(x,y)}{\nu(x_{0}, x)}=\frac{\nu(x_{0}, y)}{\nu(x_{0},x)},$$
and $$\aligned\frac{\nu^{2}(t, y)}{\nu(x_{0}, y)}\geqslant\frac{\nu^{2}(t, x)}{\nu(x_{0},x)}\endaligned,$$
So the inequality \eqref{f1.1.8} follows, which implies that
 $H_2$ is non-decreasing in $x$.

Now, we compute $\bar{\delta}_{1}$. Noticing that
$$\begin{aligned}{\big\|f_{1}^{x_{0},x_{1}}\big\|}^{2}&=\int_{0}^{x_{0}}\nu^{2}(x_{0},x_{1})\mu(\d x)+\int_{x_{0}}^{x_{1}}\nu^{2}( x, x_{1})\mu(\d x)\\
&=\nu^{2}(x_{0},x_{1})\mu(0,x_{0})+\int_{x_{0}}^{x_{1}}\nu^{2}( x, x_{1})\mu(\d x),\\
D(f_{1}^{x_{0},x_{1}})&=\int_{x_{0}}^{x_{1}}e^{C(s)}\big({(f_{1}^{x_{0}, x_{1}}})'(s)\big)^{2}\d s=
\int_{x_{0}}^{x_{1}}e^{-C(s)}\d s=\nu(x_{0}, x_{1}). \end{aligned}$$
we obtain
 $$\aligned\bar{\delta}_{1}&=\sup\limits_{x_{0},x_{1}:x_{0}<x_{1}} \frac{{\|f_{1}^{x_{0},x_{1}}\|}^{2}}{D(f_{1}^{x_{0}, x_{1}})}
\\&=\sup\limits_{x_{0},x_{1}:x_{0}<x_{1}}\bigg(\nu(x_{0}, x_{1})\mu(0, x_{0})+\frac{1}{\nu(x_{0},x_{1})}\int_{x_{0}}^{x_{1}}\!\nu^{2}(x,x_{1})\mu(\d x)\bigg). \endaligned$$
Comparing this with the expression of $\delta_{1}'$ in \eqref{f1.1.7}, we obtain $\bar{\delta}_{1}=\delta_{1}'$.
\medskip

(f) At last, we show that $\delta_{1}'\leqslant 2\delta$. Without loss of generality, assume that $\delta<\infty$.
Using the integration by parts formula, we have
$$\begin{aligned}
\int_{0}^{x}\nu^{2}(s\vee x_{0}, D)\mu(\d s)
&=\nu^{2}(x,D)\mu(0, x)+2\int_{x_{0}}^{x}\mu(0, s)\nu(s, D)e^{-C(s)}\d s\\
&\leqslant \delta\nu(x, D)+2\delta\int_{x_{0}}^{x}e^{-C(s)}\d s,\qquad x\geqslant x_{0}.
\end{aligned}$$
 By letting $x\rightarrow D$, we obtain $$\int_{0}^{D}\nu^{2}(s\vee x_{0}, D)\mu(\d s)\leqslant2\delta\nu(x_{0},D),$$ or equivalently,
$$\aligned
\nu(x_{0},D)\mu(0, x_{0})+\frac{1}{\nu(x_{0}, D)}\int_{x_{0}}^{D}\nu^{2}(s, D)\mu(\d s)
\leqslant 2\delta,\qquad x_{0}\in(0,D).\endaligned$$
Making supremum with respect to $x_{0}\in(0, D)$ on the both sides of the inequality, the assertion that $\delta_{1}'\leqslant 2\delta$ follows from \eqref{f1.1.17} immediately.$\qquad\Box$

\newpage
\subsection*{Appendix B\quad Complement of the proofs in section 4}
\subsection*{B.1\quad Proof of theorem $\ref{1t1.1}$}

Similar to the ND situation,  we adopt two circle arguments follows.
 \begin{eqnarray}
\tilde{\lambda}_{0}\hskip-0.6cm&&\geqslant\lambda_{0}\label{cr1}\\
&&\geqslant\sup\limits_{f\in\mathscr{F}_{I\!I}}\inf\limits_{x\in(0, D)}I\!I(f)(x)^{-1}\!\!=\sup\limits_{f\in\mathscr{F}_{I}}\inf\limits_{x\in(0, D)}I\!I(f)(x)^{-1}\!\!\nonumber\\&&\hskip1.5cm=\sup\limits_{f\in\mathscr{F}_{I}}\inf\limits_{x\in(0,D)}I(f)(x)^{-1}\hskip 0.8cm\label{cr2}\\
&&\geqslant\sup\limits_{h\in\mathscr{H}}\inf\limits_{x\in(0, D)}R(h)(x)\label{cr3}\\
&&\geqslant\tilde{\lambda}_{0}\label{cr4}. \end{eqnarray}
and
\begin{eqnarray}
\lambda_{0}\hskip-0.6cm&&\leqslant\inf\limits_{f\in\mathscr{\widetilde{F}}_{I\!I}\bigcup\mathscr{\widetilde{F}}_{I\!I}'}
\sup\limits_{x\in (0,D)}I\!I(f)(x)^{-1}\label{cr5}\\
&&\leqslant\inf\limits_{f\in\mathscr{\widetilde{F}}_{I\!I}}\sup\limits_{x\in
(0,D)}I\!I(f)(x)^{-1}=\!\! \inf\limits_{f\in\mathscr{\widetilde{F}}_{I}}\sup\limits_{x\in
(0,D)}I\!I(f)(x)^{-1}\nonumber\\&&\hskip1.5cm=
\inf\limits_{f\in\mathscr{\widetilde{F}}_{I}}\sup\limits_{x\in(0, D)}I(f)(x)^{-1}\hskip 0.8cm\label{cr6}\\
&&\leqslant \inf\limits_{h\in\mathscr{\widetilde{H}}}\sup\limits_{x\in(0,
D)}R(h)(x)\label{cr7}\\
&&\leqslant\lambda_{0}\label{cr8}.
\end{eqnarray}
The assertions below are proved in \rf{r5}{Theorem 1.1} and \rf{r6}{Chapter 6}.
\begin{equation}\label{cr9}
\lambda_{0}\geqslant\sup\limits_{f\in\mathscr{F}_{I\!I}}\inf\limits_{x\in(0, D)}I\!I(f)(x)^{-1}=\sup\limits_{f\in\mathscr{F}_{I}}\inf\limits_{x\in(0,D)}I(f)(x)^{-1}=\sup\limits_{f\in\mathscr{F}_{I}}\inf\limits_{x\in(0, D)}I\!I(f)(x)^{-1};
\end{equation}
\begin{equation}
\lambda_{0}\leqslant\inf\limits_{f\in\mathscr{\widetilde{F}}_{I\!I}}\sup\limits_{x\in
(0,D)}I\!I(f)(x)^{-1}=
\inf\limits_{f\in\mathscr{\widetilde{F}}_{I}}\sup\limits_{x\in(0, D)}I(f)(x)^{-1}= \inf\limits_{f\in\mathscr{\widetilde{F}}_{I}}\sup\limits_{x\in
(0,D)}I\!I(f)(x)^{-1}.\label{cr10}
\end{equation}
Actually, from \cite{r5,r6}, it is known that $$\aligned
\lambda_{0}\geqslant\sup\limits_{f\in\mathscr{F}_{I\!I}}\inf\limits_{x\in(0, D)}I\!I(f)(x)^{-1}=\sup\limits_{f\in\mathscr{F}_{I}}\inf\limits_{x\in(0,D)}I(f)(x)^{-1},\endaligned$$  and $$\aligned\lambda_{0}\leqslant\inf\limits_{f\in\mathscr{\widetilde{F}}_{I\!I}}\sup\limits_{x\in
(0,D)}I\!I(f)(x)^{-1}= \inf_{f\in\mathscr{\widetilde{F}}_{I}}\sup_{x\in
(0,D)}I(f)(x)^{-1}.\endaligned$$
Thus, \eqref{cr9}  holds since $\mathscr{F}_{I}\subseteq\mathscr{F}_{I\!I}$ and $\sup_{x\in(0, D)}I\!I(f)(x)\leqslant\sup_{x\in(0, D)}I(f)(x)$ by Cauchy's mean value theorem.
 \eqref{cr10} holds for the similar reason: $$\aligned\widetilde{\mathscr{F}}_{I}\subseteq\widetilde{\mathscr{F}}_{I\!I}\qquad \text{ and }\qquad\inf_{x\in(0, D)}I\!I(f)(x)\geqslant\inf_{x\in(0, D)}I(f)(x).\endaligned$$
 In particular, we have known  \eqref{cr1}, \eqref{cr2} and \eqref{cr6}  since the inequalities in \eqref{cr1} and \eqref{cr6} are obvious.
It remains to prove \eqref{cr3}--\eqref{cr5}, \eqref{cr7} and \eqref{cr8}.
 \medskip

We now begin to work on the additional part of the proof under the assumption that $a, b\in{\scr C}$ except proof (c) below.
\medskip

  (a) Prove that $\sup_{f\in\mathscr{F}_{I\!I}}\inf_{x\in(0, D)}I\!I(f)(x)^{-1} \geqslant
\sup_{h\in\mathscr{H}}\inf_{x\in(0, D)}R(h)(x)$.
\medskip

 Given $h\in\mathscr{H},$ let $g(x)=g(\varepsilon)\exp\big[\int_{\varepsilon}^{x}h(u)\d u\big], x\in(0,D)$ for a fixed $\varepsilon>0$. Then $g\in{\scr C}^{2}(0,D)\cap {\scr C}[0,D]$, $$\aligned g(0)=0,\qquad g'>0,\qquad\text{and}\quad h=\frac{g'}{g}\quad\text{on}\;(0,D).\endaligned$$ Furthermore,
$$\aligned R(h)=-(ah^{2}+bh+ah')=-\frac{L g}{g}. \endaligned$$

To show that
$$\aligned\inf\limits_{x\in(0,D)}R(h)(x)\leqslant\sup\limits_{f\in\mathscr{F}_{I\!I}}\inf\limits_{x\in(0,D)}I\!I(f)(x)^{-1}\qquad\text{for every}\quad h\in\mathscr{H},\endaligned$$
without loss of generality, assume that $\inf_{x\in(0, D)}R(h)(x)>0$. This ensures $$\aligned f:=-(ag''+bg')=gR(h)>0.\endaligned$$ Then  $f>0$ and
$f\in{\scr C}[0, D]$ since $a, b\in\scr{C}[0,D]$.  Since $f\!=\!-(ag''+bg')\!=\!Lg$ and $g'(D)\geqslant0$, by \eqref{fnew1}, we have
$$g'(s)\geqslant g'(s)-g'(D)= e^{-C(s)}\int_{s}^{D}f\d\mu,\qquad s\in(0,D).$$
Moreover,  we obtain
$$g(x)\geqslant\int_{0}^{x}\nu(\d s)\int_{s}^{D}f\d\mu=f(x)I\!I(f)(x),\qquad x\in(0,D)$$
since $g(0)=0$.
Thus, $$\aligned R(h)(x)^{-1}\geqslant\frac{g(x)}{f(x)}\geqslant I\!I(f)(x),\qquad x\in(0,D).\endaligned$$
Therefore,
$$\inf\limits_{x\in(0, D)}R(h)(x)\leqslant\inf\limits_{x\in(0, D)}I\!I(f)(x)^{-1}
\leqslant\sup\limits_{f\in\mathscr{F}_{I\!I}}\inf\limits_{x\in(0,D)}I\!I(f)(x)^{-1}.$$ The assertion follows since $h$ is arbitrary.
\medskip

(b) Prove that
$\sup_{h\in\mathscr{H}}\inf_{x\in(0,D)}R(h)(x)\geqslant\tilde{\lambda}_{0}$.
\medskip

Firstly, we show that $\sup_{h\in\mathscr{H}}\inf_{x\in(0,D)}R(h)(x)\geqslant0$. For a given positive $f\in L^{1}(\mu),$ let $g=fI\!I(f).$ Then $$\aligned g'(x)=e^{-C(x)}\int_{x}^{D}f\d\mu>0.\endaligned$$
 Let $h={g'}/{g}.$  By simple calculation, we get $$\aligned -f=ag''+bg'\quad\text{and}\quad R(h)=-(ah^{2}+bh+ah')=-\frac{Lg}{g}=\frac{f}{g}>0\;\;\text{ on } (0,D).\endaligned$$  This implies $\inf_{x\in(0,D)}R(h)(x)\geqslant0$ and the required assertion follows.

When $\lambda_{0}>0,$ It was proved in \rf{r4}{Theorem 2.2} and \rf{r6}{Proof (d) of Theorem 3.7}(also mentioned in the proofs of \rf{r5}{Theorem 1.2}) that the eigenfunction of $\lambda_{0}$ is strictly increasing.  Even though $\tilde{\lambda}_{0}$ could be formally bigger than $\lambda_{0}$, the same proofs still work for the eigenfunction $g$ of $\tilde{\lambda}_{0}$ since the constructed function $g$ used there satisfies $g=g(\cdot\wedge x_{0})$ for some $x_{0}\in(0,D).$ Hence, there exists an eigenfunction $g$ such that $$\aligned Lg=-\tilde{\lambda}_{0}g,\quad g(0)=0,\quad\text{and}\quad g\in{\scr C}^{2}(0,D)\cap {\scr C}[0,D].\endaligned$$ Let $h={g'}/{g}$. Then $h\in{\scr C}^{1}(0, D)\cap {\scr C}[0,D]$, $h\in\mathscr{H}$ and
 $$\aligned R(h)(x)=-\frac{Lg(x)}{g(x)}=\tilde{\lambda}\qquad\text{for}\quad x\in(0, D).\endaligned$$ So the assertion follows.
\medskip

(c) Prove that $\lambda_{0}\leqslant\inf_{ f\in\widetilde{\mathscr{F}}_{I\!I}\cup\widetilde{\mathscr{F}}'_{I\!I}}\sup_{x\in(0,D)}I\!I(f)(x)^{-1}.$
\medskip

 When $D=\infty$, this is almost done in the original proof of \rf{r5}{Theorem1.1} except that one requires an additional condition $g\in L^{2}(\mu),$ provided $x_{0}=\infty$ is allowed. This is the reason why the set $\widetilde{\mathscr{F}}'_{I\!I}$ is added. Anyhow the proof is similar to that of Theorem \ref{th1} presented in Section 3.

\medskip
(d) Prove that $\inf_{f\in\widetilde{\mathscr{F}}_{I\!I}}\sup_{x\in(0,D)}I\!I(f)(x)^{-1}\leqslant\inf_{h\in\widetilde{\mathscr{H}}}\sup_{x\in(0,D)}\mathscr{R}(h)(x).$
\medskip

Firstly, for $h\in\widetilde{\mathscr{H}},$  $\exists x_{0}\in (0,D)$ such that $h|_{(0,x_{0})}>0$, \mbox{$\int_{0+}h(u)\d u=\infty$}, and $h|_{[x_{0},D]}=0$. Similar to the proof (b) above,  given $g$, we change the form of $R(h)$ on $(0,x_{0})$. Thus,
$$R(h)=\begin{cases}-\big({ag''+bg'}\big)\big/{g},&\text{on } (0,x_{0});\\
0,&\text{otherwise}.\end{cases}$$

Next, for $h\in\widetilde{\mathscr{H}},$ let
$$\aligned f(x)=[gR(h)](x)=-ae^{-C}\big(e^{C}g'\big)'(x)\;\;\text{for}\;\; x\leqslant x_{0};\quad f(x)=f(x_{0})\;\;\text{for}\;\; x>x_{0}.\endaligned$$
Then $f\in\widetilde{\mathscr{F}}_{I\!I}$ since $a, b\in{\scr C}[0,D]$, and
$$
\aligned
e^{C(x)}g'(x)=\int_{x}^{x_{0}}f\d\mu+e^{C(x_{0})}g'(x_{0})\quad\text{for }x\leqslant x_{0};\qquad e^{C(x_{0})}g'(x_{0})=\int_{x_{0}}^{D}f\d\mu.
\endaligned$$
Moreover,  we have $g'(x)= e^{-C(x)}\int_{x}^{D}f\d\mu$,  which implies that $$\aligned g(x)=\int_{0}^{x}\nu(\d s)\int_{s}^{D}f\d\mu\quad\text{and}\quad R(h)(x)^{-1}=\frac{g(x)}{f(x)}\leqslant I\!I(f)(x)\;\;\text{for}\; x\in(0,x_{0}).\endaligned$$
Therefore, we get
 $$\aligned
 \sup_{x\in(0,x_{0})}R(h)(x)&\geqslant\sup_{x\in(0,x_{0})} I\!I(f)(x)^{-1}\\&\geqslant\inf_{f\in\widetilde{\mathscr{F}_{I\!I}},f=f_{\cdot\wedge x_{0}}}\sup_{x\in(0,x_{0})} I\!I(f)(x)^{-1}\\&\geqslant\inf_{f\in\widetilde{\mathscr{F}_{I\!I}}}\sup_{x\in(0,D)} I\!I(f)(x)^{-1}.
\endaligned$$
Furthermore, we obtain
$$\aligned
\inf\limits_{f\in\widetilde{\mathscr{F}}_{I\!I}}\sup\limits_{x\in(0,D)}I\!I(f)(x)^{-1}\leqslant\inf\limits_{h\in\widetilde{\mathscr{H}}}\sup\limits_{x\in(0,D)}R(h)(x).
\endaligned$$
\medskip

(e) Prove that
$\inf\limits_{h\in\widetilde{\mathscr{H}}}\sup\limits_{x\in(0,D)}{R}(h)(x)\leqslant\lambda_{0}.$
\medskip

Recall the definition of $\lambda_{0}:$
$$\aligned\lambda_{0}=\inf\big\{D(f):\,&\mu\big(f^{2}\big)=1,\, f(0)=0,\,f=f(\cdot\wedge x_{0}), \,f\in {\scr C}^1(0,x_{0})\cap {\scr C}[0,x_{0}]\\&\text{ for some } x_{0}\in (0,D)\big\}=:\tilde{\lambda}_{0}.\endaligned$$
Let $p_{n}\uparrow D$ and denote by $\lambda_{0}^{(0, p_{n})}$ the corresponding eigenvalue determined by $L|_{(0, p_{n})}$ (The same as the proof of \rf{r5}{Theorem 1.2}. Then $\lambda_{0}^{(0, p_{n})}\downarrow\lambda_{0}$ by using the proof of \rf{r1}{Lemma 5.1}.$\qquad \Box$

\medskip
\subsection*{B.2\quad  Proof of Theorem $\ref{1t1.3}$}
 (a) We remark that the sequence $\big\{f_{n}^{(x_{0})}\big\}_{n\in\mathbb{N}}$ is clearly contained in $\widetilde{\mathscr{F}}_{I}$.  But the modified sequence used in \rf{r5}{Theorem 1.2}:
$$
\aligned
\widetilde{f}_{1}^{(x_{0})}=\varphi(\cdot\wedge x_{0}),\qquad \widetilde{f}_{n}^{(x_{0})}=\widetilde{f}_{n-1}^{(x_{0})}(\cdot \wedge x_{0})I\!I(\widetilde{f}_{n-1}^{(x_{0})}(\cdot\wedge x_{0})),\ \ n\geqslant2,
\endaligned$$
is usually not contained in $\widetilde{\mathscr{F}}_{I\!I}$. However,
$$\aligned
\delta_{n}'&=\sup\limits_{x_{0}\in(0,D)}\inf\limits_{x\in(0,D)}I\!I\big(f_{n}^{(x_{0})}\big)(x)\\&=\sup\limits_{x_{0}\in(0,D)}\inf\limits_{x\in(0,x_{0})}I\!I\big(f_{n}^{(x_{0})}\big)(x)\\
&=\sup\limits_{x_{0}\in(0,D)}\inf\limits_{x\in(0,x_{0})}I\!I\big(\widetilde{f}_{n}^{(x_{0})}(\cdot\wedge x_{0})\big)(x)\\&=\sup\limits_{x_{0}\in(0,D)}\inf\limits_{x\in(0,D)}I\!I\big(\widetilde{f}_{n}^{(x_{0})}(\cdot\wedge x_{0})\big)(x).
\endaligned$$
Here in the last step we have used the convention that $1/0=\infty.$ Hence, these two sequences produce the same $\{\delta_{n}'\}$.
 The assertions about $\delta_{n}$ and $\delta_{n}'$ were proved in \rf{r5}{Theorem 1.2}.
\medskip

(b) Prove that $\bar{\delta}_{n+1}\geqslant\delta'_{n}$ and $\bar{\delta}_{n}^{-1}\geqslant\lambda_{0}\,(n\geqslant1)$.

The assertion of $\bar{\delta}_{n}^{-1}\geqslant\lambda_{0}$ is obvious since every function in $\big\{f_{n}^{(x_{0})}:n\geqslant1\big\}$ is a  test function of $\tilde{\lambda}_{0}$ and $\tilde{\lambda}_{0}=\lambda_{0}$.
Similar to the ND case, it is easy to see that $\bar{\delta}_{n+1}\geqslant\delta'_{n}$, which is a consequence of  the proof  of \rf{r5}{Theorem 1.1}.

Indeed, when  proving $\lambda_{0}\leqslant\inf_{f\in\mathscr{F}_{I\!I}}\sup_{x\in(0,D)}I\!I(f)(x)^{-1}$ \big(i.e $\xi_{0}'$ there\big), we know that $g=[fI\!I(f)](\cdot\wedge x_{0})$
satisfies ${D(g)}/{\mu(g^{2})}\leqslant \sup_{x\in(0,D)}I\!I(f)(x)^{-1}.$ By the relation between  $f_{n}^{(x_{0})}$ and $f_{n+1}^{(x_{0})}$,
we have $$\frac{\big\|f_{n+1}^{(x_{0})}\big\|}{D(f_{n+1}^{(x_{0})})}\geqslant\inf\limits_{x\in(0,D)}I\!I\big(f_{n}^{(x_{0})}\big)(x),$$ which implies that $$\aligned\bar{\delta}_{n+1}=\sup\limits_{x_{0}\in(0,D)}\frac{\big\|f_{n+1}^{(x_{0})}\big\|}{D\big(f_{n+1}^{(x_{0})}\big)}\geqslant\delta'_{n}. \endaligned$$
The assertion that $\delta_{1}'=\bar{\delta}_{1}$ is proved in the appendix B.3 below.$\qquad\Box$
\medskip

\subsection*{B.3\quad Proof of Corollary $\ref{1c1.1}$}

The degenerated case that $\mu(0,D)=\infty$ is trivial since $\lambda_{0}=0$  and $$\delta=\delta_{1}=\delta_{1}'=\infty.$$
The main assertion of  Corollary \ref{1c1.1}  is a consequence of Theorem \ref{1t1.3}. Here, we compute  $\delta_{1}$, $\delta_{1}'$ and prove that $\delta_{1}'\in[\delta, 2\delta]$.
Compute $\delta_{1}$ first.

Since
$$
\aligned
\int_{0}^{x}\nu(\d t)\int_{t}^{D}\sqrt{\varphi}\,\d\mu
&=\int_{0}^{x}\nu(\d t)\int_{t}^{x}\!\!\sqrt{\varphi}\,\d\mu+\int_{0}^{x}\nu(\d t)\int_{x}^{D}\sqrt{\varphi}\,\d\mu\\
&=\int_{0}^{x}\sqrt{\varphi(s)}\,\mu(\d s)\int_{0}^{s}\d\nu+\varphi(x)\int_{x}^{D}\sqrt{\varphi}\,\d\mu\\
&=\int_{0}^{x}\sqrt{\varphi}\,\varphi\d\mu+\varphi(x)\int_{x}^{D}\sqrt{\varphi}\,\d\mu\\
&=\int_{0}^{D}\sqrt{\varphi(s)}\,\varphi(s\wedge x)\,\mu(\d s),
\endaligned$$
we have $$\aligned\delta_{1}&=\sup_{x\in(0,D)}I\!I(\sqrt{\varphi} )(x)\\&=\sup\limits_{x\in(0,D)}\bigg(\frac{1}{\sqrt{\varphi(x)}}\int_{0}^{x}\!\!\sqrt{\varphi}\,\varphi\,\d\mu
+\sqrt{\varphi(x)}\int_{x}^{D}\!\!\sqrt{\varphi}\d\mu\bigg)\\&=\sup_{x\in(0,D)}\frac{1}{\sqrt{\varphi(x)}}\int_{0}^{D}\sqrt{\varphi(s)}\,\varphi(s\wedge x)\,\mu(\d s).\endaligned$$

Now, we compute  $\delta_{1}'$. Note that
$$\aligned
I\!I(f_{1}^{(x_{0})})(x)=\frac{1}{\varphi(x\wedge x_{0})} \int_{0}^{x}e^{-C(t)}\d t\int_{t}^{D}\varphi(s\wedge x_{0})\mu(\d s).
\endaligned$$
The right-hand side is clearly increasing in $x$ for $x\geqslant x_{0}$ and decreasing for $x\leqslant x_{0}.$ Hence, $I\!I\big(f_{1}^{(x_{0})}\big)$ achieves its minimum at $x=x_{0}$. By exchanging the order of the integrals, its minimum is equal to $$
\aligned
\frac{1}{\varphi(x_{0})}\int_{0}^{D}\varphi^{2}(s\wedge x_{0})\mu(\d s).
\endaligned$$  So
$$\aligned\delta_{1}'=\sup_{x_{0}\in(0,D)}\inf_{x\in(0,D)}I\!I\big(f_{1}^{(x_{0})}\big)(x)=\sup_{x_{0}\in(0,D)}\frac{1}{\varphi(x_{0})}\int_{0}^{D}\varphi^{2}(s\wedge x_{0})\mu(\d s).\endaligned$$

Next, following the proof in the discrete case \cite{r2}, we have $$\aligned
D\big(f_{1}^{(x_{0})}\big)=\int_{0}^{D}e^{C(x)}[\varphi'(x_{0}\wedge x)]^{2}\d x=\int_{0}^{x_{0}}e^{C(x)}\big(e^{-C(x)}\big)^{2}\d x=\varphi(x_{0}),
\endaligned$$ and
$$
\aligned
\mu\big(f_{1}^{(x_{0})}\big)=\int_{0}^{D}\varphi^{2}(x_{0}\wedge x)\mu(\d x).
\endaligned$$
Thus,
$$\aligned
\bar{\delta}_{1}=\sup\limits_{x_{0}\in(0,D)}\frac{\mu\big[\big(f_{1}^{(x_{0})}\big)^{2}\big]}{D\big(f_{1}^{(x_{0})}\big)}=\delta_{1}'.
\endaligned$$

At last, we prove that $\delta_{1}'\in[\delta,2\delta]$.

Following the corresponding proof in the discrete case \rf{r2}{Corollary 4.4},
we have
$$\delta_{1}'=\sup\limits_{ x_{0}\in(0,D)}\frac{1}{\varphi(x_{0})}\int_{0}^{D}\varphi^{2}(s\wedge x_{0})\mu(\d s)\geqslant\sup\limits_{x\in(0,D)}\varphi(x_{0})\mu(x_{0},D)=\delta.$$
On the other hand,  using the integration formula by parts, for $x<x_{0}$, we have
$$\int_{x}^{x_{0}}\varphi^{2}(s)\mu(\d s)=-\varphi^{2}(s)\mu(s,D)\big|_{x}^{x_{0}}+2\int_{x}^{x_{0}}\varphi(s)\varphi'(s)\mu(s,D)\d s.$$
So
$$\aligned
\int_{x}^{D}\varphi^{2}(s\wedge x_{0})\mu(\d s)&=\int_x^{x_{0}}\varphi^{2}(s)\mu(\d s)+\varphi^{2}(x_{0})\mu(x_{0},D)\\
&=\varphi^{2}(x)\mu(x,D)+2\int_{x}^{x_{0}}\varphi(s)\varphi'(s)\mu(s,D)\d s\\
&\leqslant \delta \varphi(x)+2\delta\int_{x}^{x_{0}}e^{-C(s)}\d s\\&\rightarrow 2\delta \varphi(x_{0})\qquad \text{as}\ x\rightarrow0.
\endaligned$$
Thus ${\int_{0}^{D}\!\varphi^{2}(s\wedge x_{0})\mu(\d s)}/{\varphi(x_{0})}\!\leqslant\! 2\delta$ and
the assertion $\delta_{1}'\!\leqslant\!2\delta$ follows
immediately.$\; \Box$


\begin{enumerate}
{\footnotesize

\bibitem{r4}\label{r4} Chen M F. Analytic proof of dual variational formula for the first eigenvalue in dimension one.  Sci in China, A, 1999, 42(8): 805--815\\[-6.5mm]

\bibitem{r3}\label{r3} Chen M F. Explicit bounds of the first eigenvalue. Sci Chin Ser, A, 2000, 43(10): 1051--1059\\[-6.5mm]

\bibitem{r5}\label{r5} Chen M F. Variational formulas and approximation theorems for the first eigenvalue in dimension one.  Sci in China, A, 2001, 44(4): 409--418\\[-6.5mm]

\bibitem{r6}\label{r6} Chen M F. Eigenvalues, Inequalities, and Ergodic theory. New York: Springer, 2005\\[-6.5mm]

\bibitem{r2}\label{r2} Chen M F. Speed of stability for birth-death process. Front Math China, 2010, 5(3): 379--516\\[-6.5mm]

\bibitem{rn1}\label{rn1} Chen M F. Basic estimates and stability rate for one-dimensional diffusion. Chapter 6 in ``Probability Approximations and Beyond", Lecture notes in statistics 205, 2012, 75--99\\[-6.5mm]

\bibitem{rn2}\label{rn2} Chen M F. General estimate of the first eigenvalue on manifolds. Front Math China, 2011, 6(6): 1025--1043\\[-6.5mm]

\bibitem{rn3}\label{rn3} Chen M F. Lower bounds of the principle eigenvalue in dimension one.  Front Math China, to  appear\\[-6.5mm]

\bibitem{r1}\label{r1} Chen M F,  Wang F Y.  Estimation of spectral gap for elliptic operators. Trans Amer Math Soc, 1997, 349(3): 1239--1267. See also book [4] in\\
    http://math.bnu.edu.cn/\~{}chenmf/main$\_$eng.htm\\
    for some supplement.\\[-6.5mm]

\bibitem{wj}\label{wj} Wang J. First eigenvalue of one-dimensional diffusion processes. Elect Comm Prob, 2009, 14: 232-244\\[-6.5mm]

\bibitem{zettl}\label{zettl} Zettl A. Sturm-Liouville Theory, AMS, 2005

}\end{enumerate}
\end{document}